\documentclass[aos,MSNbibl,nameyear,seceqn,dvips]{arximspdf}
\usepackage{url,breakurl}
\doi{10.1214/14-AOS1258} % kopijuoti is laisko
\volume{42}
\issue{6}
\pubyear{2014}
\firstpage{2441}
\lastpage{2468}
\docsubty{FLA}

\makeatletter

\newtheorem{thmm}{Theorem}[section]
\newtheorem{coro}{Corollary}[section]

\newtheorem{pro}{Proposition}[section]
\newproclaim{rem}{Remark}[section]
\newproclaim{defi}{Definition}[section]
\newproclaim{exm}{Example}[section]

\makeatother

\begin{document}
\begin{frontmatter}

\title{Asymptotic theory of generalized information criterion for geostatistical regression model selection}
\runtitle{Asymptotics for geostatistical regression model selection\hspace*{4pt}}

\begin{aug}
\author[A]{\fnms{Chih-Hao}~\snm{Chang}\ead[label=e1]{jhow@stat.sinica.edu.tw}},
\author[B]{\fnms{Hsin-Cheng}~\snm{Huang}\corref{}\thanksref{T1}\ead[label=e2]{hchuang@stat.sinica.edu.tw}\ead[label=u1,url]{http://www.stat.sinica.edu.tw/}}
\and
\author[B]{\fnms{Ching-Kang}~\snm{Ing}\thanksref{T2}\ead[label=e3]{cking@stat.sinica.edu.tw}}
\runauthor{C.-H. Chang, H.-C. Huang and C.-K. Ing}
\affiliation{National University of
Kaohsiung, Academia Sinica and Academia Sinica}
\address[A]{C.-H. Chang\\
Institute of Statistics\\
National University of Kaohsiung\\
700 Kaohsiung University Road\\
Kaohsiung 811\\
Taiwan\\
\printead{e1}}
\address[B]{H.-C. Huang\\
C.-K. Ing\\
Institute of Statistical Science\\
Academia Sinica\\
128 Academia Road Section~2\\
Taipei 115\\
Taiwan\\
\printead{e2}\\
\phantom{E-mail:\ }\printead*{e3}\\
\printead{u1}}
\end{aug}
\thankstext{T1}{Supported by National Science Council of Taiwan under Grant 100-2628-M-001-004-MY3.}
\thankstext{T2}{Supported by National Science Council of Taiwan under Grant 99-2118-M-001-008-MY2 and
the Academia Sinica Investigator Award.}

% HISTORY:
\received{\smonth{3} \syear{2014}}
\revised{\smonth{7} \syear{2014}}

% ABSTRACT
\begin{abstract}
Information criteria, such as Akaike's information criterion and
Bayesian information criterion are often applied in model selection.
However, their asymptotic behaviors for selecting geostatistical
regression models have not been well studied, particularly under the
fixed domain asymptotic framework with more and more data observed
in a bounded fixed region. In this article, we study the generalized
information criterion (GIC) for selecting geostatistical regression
models under a more general mixed domain asymptotic framework.
Via uniform convergence developments of some statistics,
we establish the selection consistency and the asymptotic loss efficiency of GIC
under some regularity conditions, regardless of whether the covariance model is correctly
or wrongly
specified. We further provide specific examples
with different types of explanatory variables that satisfy the conditions. For example,
in some situations, GIC is selection consistent, even when some
spatial covariance parameters cannot be estimated consistently.
On the other hand, GIC fails to
select the true polynomial order consistently under the fixed domain
asymptotic framework. Moreover, the growth rate of the domain and
the degree of smoothness of candidate regressors in space are shown to play key
roles for model selection.
\end{abstract}

% KEYWORDS

\end{frontmatter}

%s1 #&#
\section{Introduction}\label{sec1}
With the advent of data collection technologies, more and more data,
such as remote sensing data or environmental monitoring data, are
collected in space and managed by geographical information systems.
In many applications, a response of interest is observed on a set of
sites in space, and it is of interest to apply a geostatistical
regression model to predict the response at unsampled sites with the
aid of auxiliary/explanatory variables. For example, in precision
agriculture, it is of interest to predict crop yield based on some
explanatory variables involving, for example, climatic conditions,
soil types, fertilizers, cropping practices, weeds and topographic
features. Not only do we aim to identify the important
explanatory variables, but the precision of yield also depends on
how well the explanatory variables are chosen, which if not chosen properly, may result in poor
performance, particularly when the number of explanatory variables
is large. Clearly, model selection is essential in geostatistics.

There are two different asymptotic frameworks in geostatistics. One
is called the increasing domain asymptotic framework, where the
observation region grows with the sample size. The other is called
the fixed domain asymptotic (or infill asymptotic) framework, where
the observation region is bounded and fixed with more and more data
taken more densely in the region. It is known that the two
frameworks lead to possibly different asymptotic behaviors in
covariance parameter estimation. However, little is known about
their effects on model selection. In general, asymptotic behaviors
of the estimated parameters under the increasing domain framework
are more standard. For example, the maximum likelihood estimates of
covariance parameters are typically consistent and asymptotically
normal when fitted by a correct model [\citet{Mardia1984}]. In
contrast, not all covariance parameters can be estimated
consistently under the fixed domain asymptotic framework, even for
the simple exponential covariance model in one dimension with no
consideration of explanatory variables [\citet{Ying1991}; Chen, Simpson and Ying
(\citeyear{Chen2000})]. The readers are refereed to \citet{Stein1999} for more
details regarding fixed domain asymptotics. Some discussion
concerning which asymptotic framework is more appropriate can also
be found in \citet{Zhang2005}.

Many model selection methods have been applied in geostatistical
regression, such as Akaike's information criterion [AIC, \citet{Akaike1973}],
Bayesian information criterion [BIC, Schwartz (\citeyear{Schwart1978})], the
generalized information criterion [GIC, \citet{Nishii1984}] and the cross
validation method [\citet{Stone1974}]. Note that GIC contains a range of
criteria, including both AIC and BIC, governed by a tuning
parameter. Although theoretical properties of these selection
methods have been thoroughly established in linear regression and
time series model selection [e.g., \citet{Shao1997}, McQuarrie and Tsai
(\citeyear{McQuarrie1989}), \citet{Ing2005}, \citet{Ing2007}], only limited results are
available for selecting geostatistical regression models. For
example, \citet{Hoeting2006} provided some heuristic
arguments for AIC in geostatistical model selection when the spatial
process of interest is observed with no measurement error. They show
via a simulation study that spatial dependence has to be considered,
which if ignored, may result in unsatisfactory results. \citet{Huang2007} developed a technique of estimating the mean squared
prediction error for a general spatial prediction procedure using
the generalized degrees of freedom and derived some asymptotic
efficiency results. For linear mixed models, \citet{Jiang2003}
developed some consistent procedures similar to GIC. \citet{Pu2006} derived conditions under which GIC is selection consistent.
\citet{Jiang2008} introduced a fence method for mixed
model selection and showed its consistency under some regularity
conditions. \citet{Jones2011} proposed a modified BIC, which replaces the
sample size in the penalty of the original BIC by an effective
sample size to account for correlations in linear mixed models.
\citet{Vaida2005} proposed the conditional Akaike's
information criterion (CAIC) and argued that it is more appropriate
than AIC when the focus is on subjects/clusters requiring prediction
of random effects. In addition, selection among semiparametric
regression models and penalized smoothing spline models [e.g.,
Chapter~4, Ruppert, Wand and Carroll (\citeyear{Puppert2003})] can also be formulated
in terms of random-effect selection in linear mixed models. The
asymptotic theory of AIC for this type of model was given by \citet{Shi1999}, and that for BIC was given by \citet{Bunea2004}. A recent
review of linear and generalized linear mixed model selection can
also be found in M{\"u}ller, Scealy and Welsh (\citeyear{Muller2013}).

Although the geostatistical regression model can be regarded as a
linear mixed model with one random effect, its asymptotic behavior
is surprisingly subtler than a usual linear mixed model for
the following three reasons. First, variables in a geostatistical
regression model are sampled from a spatial process, resulting in
small ``effective sample size'' unless the spatial domain is allowed
to grow quickly. Second, unlike some random-effect models with
independent random components, spatial dependence forces all
variables to depend in a complex way, making it very difficult to
handle asymptotically. Third, under the fixed domain asymptotic
framework, classical regularity conditions are generally not
satisfied, and traditional approaches for establishing asymptotic
results are typically not applicable. To the best of our knowledge,
asymptotic properties of GIC for geostatistical regression model
selection have yet to be developed, particularly under the fixed
domain asymptotic framework, where nonstandard behaviors are often
expected. In this article, we focus on GIC for geostatistical
regression model selection regardless of whether the covariance
model is correctly or wrongly specified. Although a conditional-type
criterion, such as CAIC may be more appropriate when spatial
prediction is of main interest, it is beyond the scope of this
paper. Major accomplishments are listed in the following:

\begin{longlist}[(1)]
\item[(1)] We establish a general theory of GIC for the selection consistency and
    the asymptotic loss efficiency under mild regularity conditions
    in a general mixed domain asymptotic framework, which
    includes both the fixed and increasing domain asymptotics.
    In particular, we allow the possibilities that some covariance
    parameters may converge to a nondegenerate distribution and
    the covariance model may be mis-specified.
\item[(2)] We provide some examples that satisfy
    the aforementioned regularity conditions under exponential
    covariance models in one and two dimensions, and
    demonstrate how selection consistency is affected by
    candidate regressors.
\end{longlist}

 We shall show that the asymptotic behaviors of GIC are
related to how fast the domain grows with the sample size. In
addition, some nonstandard properties of GIC under the fixed domain
asymptotic framework will be highlighted. For example, under fixed
domain asymptotics, GIC fails to identify the correct order of
polynomial consistently regardless of the tuning parameter value,
even when the underlying covariance model is correctly
specified. On the other hand, for a properly chosen tuning
parameter value, GIC is selection consistent when candidate
explanatory variables are generated from some spatial
dependent processes.

This article is organized as follows. Section~\ref{models and
criterion} gives a brief introduction of geostatistical regression
models and GIC. Our main results for the consistency and the
asymptotic loss efficiency of GIC are presented in Sections~\ref{section:selection for correct cov model} and
\ref{section:covariance model selection}. Specifically, in Section~\ref{section:selection for correct cov model}, we assume that the
covariance model is specified correctly. While in Section~\ref{section:covariance model selection}, we consider the covariance
model to be mis-specified. In Section~\ref{section:examples}, we
provide some examples that satisfy the regularity
conditions. Finally, a brief discussion is provided in Section~\ref{discussion}.

%s2 #&#
\section{Models and criteria}
\label{models and criterion}

%s2.1 #&#
\subsection{Geostatistical regression models}
\label{subsec:geo model}

Consider a spatial process, $\{S(\mathbf{s})\dvtx\break \mathbf{s}\in
D\subset\mathbb{R}^d\}$. Suppose that we observe data
$\{Z(\mathbf{s}_{n1}),\ldots,Z(\mathbf{s}_{nn})\}$ according to the following
measurement equation:
%
%e2.1 #&#
\begin{eqnarray}\label{geo data}
Z(\mathbf{s}_{ni}) &=& S(\mathbf{s}_{ni}) + \epsilon(
\mathbf{s}_{ni})
\nonumber
\\[-8pt]
\\[-8pt]
\nonumber
&=& \mu_0(\mathbf{s}_{ni})+\eta(
\mathbf{s}_{ni})+\epsilon(\mathbf{s}_{ni});\qquad i=1,\ldots,n,
\end{eqnarray}
where $\mu_0(\cdot)$ is the mean function,
$\eta(\cdot)$ is a zero-mean Gaussian spatial dependent process
with $ \sup_{\mathbf{s}\in D}\mathrm{E}(\eta^2(\mathbf{s}))<\infty$ and
$\{\epsilon(\mathbf{s}_{ni})\dvtx i=1,\ldots,n\}$ are Gaussian white-noise variables with
variance $v^2$, independent of $S(\cdot)=\mu_0(\cdot)+\eta(\cdot)$,
corresponding to measurement errors.

In addition to $Z(\mathbf{s}_{ni})$'s, we observe
$\mathbf{x}(\mathbf{s}_{ni})= (1, x_1(\mathbf{s}_{ni}),\ldots, x_{p_n}(\mathbf{s}_{ni}))'$, a
$(p_n+1)$-vector of explanatory variables, for $i=1,\ldots,n$. We
consider the geostatistical regression model
\[
Z(\mathbf{s}_{ni}) = {\mathbf x}(\mathbf{s}_{ni})'
\bolds\beta_{n}+\eta(\mathbf{s}_{ni})+ \epsilon(
\mathbf{s}_{ni});\qquad \mathbf{s}_{ni}\in D, i=1,\ldots,n,
\]
where $\bolds\beta_{n}=(\beta_0,\beta_1,\ldots,\beta_{p_n})'$.
%Here we allow the number of variables to depend on $n$ and thus is
%denoted by $p_n$. Although we also allow each $\mathbf{s}_i$ to vary
%with $n$, for simplicity, we write $\{\mathbf{s}_1,\ldots,\mathbf{s}_n\}$
%instead of $\{\mathbf{s}_{n1},\ldots,\mathbf{s}_{nn}\}$.
This model reduces
to the usual linear regression model when $\eta(\cdot)$ is absent.
Similarly to linear regression, a large model that contains many
insignificant variables may produce a large variance, resulting in
low predictive power. On the other hand, a model that ignores some
important variables may suffer from a large bias. To strike a good
balance between (squared) bias and variance, it is essential to
include only significant variables in the model. Clearly, variable
selection is essential not only in regression but also in
geostatistical regression.

We use $\alpha\subseteq\{1,\ldots,p_n\}$ to denote a model, which
consists of the indices of the corresponding explanatory variables.
Let $\mathcal{A}_n\subseteq2^{\{1,\ldots,p_n\}}$ be the set of all
candidate models with $\varnothing$ being the intercept-only model.
Let ${\mathbf X}_{n}$ be the $n\times (p_n+1)$ matrix with the $i$th row,
$\mathbf{x}(\mathbf{s}_{ni})'$; $1\leq i\leq n$. Also let ${\mathbf X}_{n}(\alpha)$ be
an $n\times (p(\alpha)+1)$ sub-matrix of ${\mathbf X}_{n}$ containing a
column $\mathbf{1}$ (corresponding to the intercept) and the columns
corresponding to $\alpha\in\mathcal{A}_n$, and $\bolds\beta_{n}(\alpha)$ be
the sub-vector of $\bolds\beta_{n}$ corresponding to ${\mathbf X}_{n}(\alpha)$. A
model $\alpha$ is said to be correct if $\mu_0({\mathbf s})$ can be
written as $ \beta_0+\sum_{j\in\alpha} \beta_j x_j({\mathbf
s})$ for all ${\mathbf s}\in D$. If there exists a correct model, we
denote the correct model having the smallest number of variables by
$\alpha_n^0= \operatorname{arg\, min}_{\alpha\in\mathcal{A}_n^0}p(\alpha)$,
where $\mathcal{A}_n^0$ is the set of all correct models.
%We allow both $\mathcal{A}_n$ and $\mathcal{A}^0_n$ to depend on $n$
%but are nested in the sense that
%$\mathcal{A}^0_k\subseteq\mathcal{A}^0_{k+1}$ and
%$\mathcal{A}_k\setminus\mathcal{A}^0_k\subseteq\mathcal{A}_{k+1}\setminus\mathcal{A}^0_{k+1}$,
%for $k=1,2,\ldots$.

The geostatistical regression model $\alpha$ can be written in a matrix form as
%
%e2.2 #&#
\begin{equation}
\qquad\mathbf{Z}_{n}=\bigl(Z(\mathbf{s}_{n1}),\ldots,Z(
\mathbf{s}_{nn})\bigr)' = {\mathbf X}_{n}(
\alpha)\bolds\beta_{n}(\alpha)+\bolds\eta_{n}+\bolds
\epsilon_{n};\qquad \alpha\in\mathcal{A}_n, \label{setup}
\end{equation}
where $\bolds\eta_{n}=(\eta(\mathbf{s}_{n1}),\ldots,\eta(\mathbf{s}_{nn}))'\sim N(0,{\bolds\Sigma}_{n\eta})$ and
$\bolds\epsilon_{n}=(\epsilon(\mathbf{s}_{n1}),\ldots,\epsilon(\mathbf{s}_{nn}))'$ $\sim N(0,v^2{\mathbf I}_{n})$
with ${\bolds\Sigma}_{n\eta}=\mathrm{E}(\bolds\eta_{n}\bolds\eta^{\prime}_{n})$
and ${\mathbf I}_{n}$ denoting the $n \times n$ identity matrix. Hence the mean and the variance
of $\mathbf{Z}_{n}$ conditional on $\mathbf{X}_{n}$ based on model
$\alpha\in\mathcal{A}_n$ are $\mathbf{X}_{n}(\alpha)\bolds\beta_{n}(\alpha)$ and
%
%e2.3 #&#
\begin{equation}
\bolds{\Sigma}_{n}(\bolds{\theta})= \bolds{\Sigma}_{n\eta}+v^2{
\mathbf I}_{n}, \label{Sigma}
\end{equation}
where $\bolds{\theta}$ is a covariance parameter
vector belonging to some parameter space $\Theta$. Throughout the
paper, we assume that $\bolds{\Sigma}_{n}(\bolds{\theta})$ is continuous on
$\bolds{\theta}\in\Theta$. We denote the true covariance matrix by
$\bolds\Sigma_{n0}$ and the true mean of $\mathbf{Z}_{n}$ conditional on $\mathbf{X}_{n}$
by $\bolds\mu_{n0}$. In other words, given $\mathbf{X}_{n}$, the data $\mathbf{Z}_{n}$ are
generated from $N(\bolds\mu_{n0},\bolds\Sigma_{n0})$.

In order to facilitate mathematical exposition,
the asymptotic results established in Sections~\ref{section:selection for correct cov
model} and \ref{section:covariance model selection}
focus only on the case where $\mathbf{X}_{n}$ is nonrandom.
These results are also valid in the almost sure sense
when $\mathbf{X}_{n}$ is random, provided that the required conditions
involving $\mathbf{X}_{n}$ hold for almost all sequences $\mathbf{X}_{n}$; $n\in\{1, 2, \ldots\}$.
We further illustrate these results in Section~\ref{section:examples}
using either random or nonrandom~$\mathbf{X}_{n}$.

%s2.2 #&#
\subsection{Generalized information criterion}

For notational simplicity, we suppress the dependence
of $\mathbf{X}_{n}, \mathbf{X}_{n}(\alpha), \bolds{\beta}_{n}, \bolds{\beta}_{n}(\alpha),
\mathbf{Z}_{n}$, $\bolds{\eta}_{n}$, $\bolds{\epsilon}_{n}$, ${\bolds\Sigma}_{n\eta}$, ${\mathbf I}_{n}$,
${\bolds\Sigma}_{n}(\bolds{\theta})$, ${\bolds\Sigma}_{n0}$, $\bolds\mu_{n0}$ and $\mathbf{s}_{ni}$
on $n$ in the rest of this paper.
To estimate $\bolds{\beta}$ and $\bolds{\theta}$, we\break consider maximum
likelihood (ML). We assume that $\bolds\Sigma^{-1}(\bolds\theta)$ and\break
$(\mathbf{X}'\bolds\Sigma^{-1}(\bolds\theta)\mathbf{X})^{-1}$ exist for
$\bolds\theta\in\Theta$. The ML estimate of $\bolds{\theta}$ based on
$\alpha\in\mathcal{A}_n$, denoted by $\hat{\bolds\theta}(\alpha)$, is
obtained by maximizing the following profile log-likelihood
function:
\begin{eqnarray*}
\ell(\alpha;\bolds\theta)& = & -\tfrac{1}{2}n\log(2\pi)-\tfrac{1}{2}
\log\det\bigl({\bolds\Sigma}({\bolds\theta})\bigr)
\\
& &{}-\tfrac{1}{2}\bigl({\mathbf Z}-\hat{\bolds\mu}(\alpha;\bolds\theta)
\bigr)'{\bolds\Sigma}^{-1} ({\bolds\theta}) \bigl({\mathbf
Z}-\hat{\bolds\mu}(\alpha;\bolds\theta)\bigr),
\end{eqnarray*}
where
$\hat{\bolds\mu}(\alpha;\bolds\theta)=\mathbf{X}(\alpha)\hat{\bolds\beta}(\alpha;\bolds\theta)$, and
\[
\hat{\bolds\beta}(\alpha;\bolds\theta) =\bigl({\mathbf X}(
\alpha)'{\bolds\Sigma}^{-1}(\bolds\theta){\mathbf X}(
\alpha)\bigr)^{-1} {\mathbf X}(\alpha)'{\bolds
\Sigma}^{-1}(\bolds\theta){\mathbf Z}.
\]
Specifically,
$\ell (\alpha;\hat{\bolds\theta}(\alpha) )=
 \sup_{\bolds\theta\in\Theta} \ell(\alpha;\bolds\theta)$, and
$\hat{\bolds{\beta}} (\alpha;\hat{\bolds{\theta}}(\alpha) )$ is the
ML estimate of $\bolds\beta(\alpha)$. For $\alpha\in\mathcal{A}_n$ and
$\bolds{\theta}\in\Theta$, let
%
%e2.4 #&#
%e2.5 #&#
\begin{eqnarray}
{\mathbf M}(\alpha;\bolds\theta)& =& {\mathbf X}(\alpha) \bigl({\mathbf X}(
\alpha)'{\bolds\Sigma}^{-1}(\bolds\theta){\mathbf X}(
\alpha)\bigr)^{-1}{\mathbf X}(\alpha)'{\bolds
\Sigma}^{-1}(\bolds\theta), \label{fn:M}
\\
{\mathbf A}(\alpha;\bolds\theta) &=& {\mathbf I}-{\mathbf M}(\alpha;\bolds
\theta). \label{fn:A}
\end{eqnarray}
Then $\hat{\bolds\mu}(\alpha;\bolds\theta)={\mathbf M}(\alpha;\bolds\theta)\mathbf{Z}$ and
$\mathbf{Z}-\hat{\bolds\mu}(\alpha;\bolds\theta)={\mathbf A}(\alpha;\bolds\theta)\mathbf{Z}$. Note that ${\mathbf M}^2(\alpha;\bolds\theta)
={\mathbf M}(\alpha;\bolds\theta)$, ${\mathbf M}(\alpha;\bolds\theta)\mathbf{X}(\alpha)=\mathbf{X}(\alpha)$, and
\begin{eqnarray*}
{\mathbf M}(\alpha;\bolds\theta)'\bolds{\Sigma}^{-1}(
\bolds{\theta}){\mathbf M}(\alpha;\bolds\theta) &=& \bolds{\Sigma}^{-1}(
\bolds{\theta}){\mathbf M}(\alpha;\bolds\theta),
\\
{\mathbf A}(\alpha;\bolds\theta)'\bolds{\Sigma}^{-1}(
\bolds{\theta}){\mathbf A}(\alpha;\bolds\theta)& =& \bolds{\Sigma}^{-1}(
\bolds{\theta}){\mathbf A}(\alpha;\bolds\theta).
\end{eqnarray*}
Therefore, by (\ref{fn:M}) and (\ref{fn:A}), the profile log-likelihood function can
also be written as
%
%e2.6 #&#
\begin{eqnarray}
\label{log like fun alpha} \ell(\alpha;\bolds\theta) &=& -\tfrac{1}{2}n\log(2\pi) -
\tfrac{1}{2}\log\det\bigl(\bolds\Sigma(\bolds\theta)\bigr) -\tfrac{1}{2}
\bolds\mu_0'\bolds\Sigma^{-1}(\bolds\theta)
\mathbf{A}(\alpha;\bolds\theta)\bolds\mu_0
\nonumber
\\
&&{} -\bolds\mu_0'\bolds\Sigma^{-1}(\bolds
\theta)\mathbf{A}(\alpha;\bolds\theta) (\bolds\eta+\bolds\epsilon) -
\tfrac{1}{2}(\bolds\eta+\bolds\epsilon)'\bolds
\Sigma^{-1}(\bolds\theta) (\bolds\eta+\bolds\epsilon)
\\
&&{} +\tfrac{1}{2}(\bolds\eta+\bolds\epsilon)'\bolds
\Sigma^{-1}(\bolds\theta)\mathbf{M}(\alpha;\bolds\theta) (\bolds\eta+
\bolds\epsilon); \qquad \alpha\in\mathcal{A}_n, \bolds{\theta}\in\Theta.\nonumber
\end{eqnarray}

To identify the smallest correct model $\alpha_n^0$, one may adopt the
GIC of \citet{Nishii1984},
%
%e2.7 #&#
\begin{equation}
\Gamma_{\tau_n}(\alpha) = -2\ell\bigl(\alpha;\hat{\bolds\theta}(\alpha)
\bigr) + {\tau_n} p(\alpha);\qquad \alpha\in\mathcal{A}_n,
\label{unknown GIC}
\end{equation}
where ${\tau_n}$ is a tuning parameter controlling the
trade-off between goodness-of-fit and the model parsimoniousness.
The criterion includes AIC (when ${\tau_n}=2$) and BIC [when
${\tau_n}=\log(n)$] as special cases, and has been widely used in
many statistical areas. The model selected by GIC based on
${\tau_n}$ is denoted by $ \hat{\alpha}_{\tau_n}
=\operatorname{arg\, min}_{\alpha\in\mathcal{A}_n}
\Gamma_{\tau_n}(\alpha)$. In the next section, we shall first
investigate GIC for variable selection when the covariance model is
correctly specified.

%s3 #&#
\section{Variable selection under a correct covariance model}
\label{section:selection for correct cov model}

The asymptotic properties of GIC will be derived in terms of the
Kullback--Leibler (KL) loss, which for $\alpha\in\mathcal{A}_n$ and
$\bolds\theta\in\Theta$ is given by
\begin{eqnarray*}
L(\alpha;\bolds\theta) &=& \int_{\mathbf{Y}\in\mathbb{R}^n}f(\mathbf{Y};\bolds
\mu_0,\bolds\Sigma_0) \log\frac{f(\mathbf{Y};\bolds\mu_0,\bolds\Sigma_0)}    {
f(\mathbf{Y};\hat{\bolds\mu}(\alpha;\bolds\theta),\bolds{\Sigma}(\bolds\theta))}\,d\mathbf{Y}
\\
&=& \frac{1}{2}\log\det\bigl(\bolds\Sigma(\bolds\theta)\bigr)-
\frac{1}{2}\log\det(\bolds\Sigma_0) +\frac{1}{2}
\operatorname{tr}\bigl(\bolds\Sigma_0\bolds\Sigma^{-1}(
\bolds\theta)\bigr)
\\
&&{} -\frac{n}{2}+\frac{1}{2}\bigl(\hat{\bolds\mu}(\alpha;\bolds
\theta)-\bolds\mu_0\bigr)'{\bolds\Sigma}^{-1}(
\bolds\theta) \bigl(\hat{\bolds\mu}(\alpha;\bolds\theta)-\bolds\mu_0
\bigr),
\end{eqnarray*}
where
$\hat{\bolds\mu}(\alpha;\bolds\theta)=\mathbf{X}(\alpha)\hat{\bolds\beta}(\alpha;\bolds\theta)$
and $f(\cdot;\bolds\mu,\bolds\Sigma)$ is the Gaussian density function
with mean $\bolds\mu$ and covariance matrix $\bolds\Sigma$. Note
that $L(\alpha;\bolds\theta)\geq 0$, for any $\alpha\in\mathcal{A}_n$ and
$\bolds\theta\in\Theta$.
When $\bolds{\mu}_0$ is known, the KL loss for $\bolds{\theta}\in\Theta$ is
given by
\[
L_0(\bolds\theta)= \tfrac{1}{2} \bigl\{\log\det\bigl(\bolds
\Sigma(\bolds\theta)\bigr)- \log\det(\bolds\Sigma_0)+\operatorname{tr}
\bigl(\bolds\Sigma_0\bolds\Sigma^{-1}(\bolds\theta)\bigr)-n
\bigr\}. %\label{fn:min KL}
\]
Then the optimal vector of $\bolds{\theta}\in\Theta$, which minimizes
the KL loss, is given by
\[
\bolds\theta_0= \mathop{\operatorname{arg\, inf}}_{\bolds\theta\in\Theta}
L_0(\bolds\theta). %\label{def:theta 0}
\]

Clearly, $\bolds{\Sigma}_0=\bolds{\Sigma}(\bolds{\theta}_0)$ and
$L_0(\bolds{\theta}_0)=0$, if the covariance model class contains the
correct model. In this case, $\bolds{\theta}_0$ is the true covariance
parameter vector of $\bolds{\theta}$. Let
$R(\alpha;\bolds\theta)=\mathrm{E}(L(\alpha;\bolds\theta))$. By
(\ref{fn:M}) and (\ref{fn:A}), we have
%
%e3.1 #&#
%e3.2 #&#
\begin{eqnarray} \label{fn:KL loss}
L(\alpha;\bolds\theta) &=& L_0(\bolds\theta) +\tfrac{1}{2}\bolds
\mu_0'\bolds\Sigma^{-1}(\bolds\theta)
\mathbf{A}(\alpha;\bolds\theta)\bolds\mu_0
\nonumber
\\[-8pt]
\\[-8pt]
\nonumber
&&{} +\tfrac{1}{2}(\bolds\eta+\bolds\epsilon)'\bolds
\Sigma^{-1}(\bolds\theta)\mathbf{M}(\alpha;\bolds\theta) (\bolds\eta+
\bolds\epsilon),
\\
\label{fn:KL risk}
R(\alpha;\bolds\theta) &=& L_0(\bolds\theta)+\tfrac{1}{2}\bolds
\mu_0'\bolds\Sigma^{-1}(\bolds\theta)
\mathbf{A}(\alpha;\bolds\theta)\bolds\mu_0
\nonumber
\\[-8pt]
\\[-8pt]
\nonumber
&&{} +\tfrac{1}{2}\operatorname{tr}\bigl( \bolds\Sigma^{-1}(\bolds
\theta)\mathbf{M}(\alpha;\bolds\theta)\bolds\Sigma_0\bigr),
\end{eqnarray}
for $\alpha\in\mathcal{A}_n$ and $\bolds\theta\in\Theta$, where
$\bolds\mu_0'\bolds\Sigma^{-1}(\bolds\theta)
\mathbf{A}(\alpha;\bolds\theta)\bolds\mu_0= \|\bolds{\Sigma}^{-1/2}(\bolds{\theta})
\mathbf{A}(\alpha;\bolds{\theta})\bolds{\mu}_0\|^2$, which results\vspace*{1pt} from using a wrong regression model, and is equal to $0$
when $\alpha\in\mathcal{A}_n^0$. In particular, for $\alpha\in\mathcal{A}_n^0$ and $\bolds\Sigma_0=\bolds\Sigma(\bolds\theta_0)$,
%
%e3.3 #&#
%e3.4 #&#
\begin{eqnarray}
L(\alpha;\bolds\theta_0) &=& \tfrac{1}{2}(\bolds\eta+\bolds
\epsilon)'\bolds\Sigma^{-1}(\bolds\theta_0)
\mathbf{M}(\alpha;\bolds\theta_0) (\bolds\eta+\bolds\epsilon),
\label{fn:KL loss 0}
\\
R(\alpha;\bolds\theta_0) &=& \tfrac{1}{2}p(\alpha).
\label{fn:KL risk 0}
\end{eqnarray}

Consider a model selection procedure $\hat{\alpha}$ that maps data to
$\alpha\in\mathcal{A}_n$. We say that $\hat{\alpha}$ is consistent if
$ \lim_{n\rightarrow\infty}P \{\hat{\alpha}=\alpha_n^0 \}=1$,
and $\hat\alpha$ is asymptotically loss efficient if
%
%e3.5 #&#
\begin{equation}
\frac{L(\hat\alpha;\hat{\bolds\theta}(\hat\alpha))}{ \min_{\alpha\in\mathcal{A}_n}
    L(\alpha;\hat{\bolds\theta}(\alpha))}  \mathop{\rightarrow}\limits^{P}1, \label{loss efficiency}
\end{equation}
as $n\rightarrow\infty$. When $\eta(\cdot)$ is absent,
geostatistical regression reduces to the usual linear regression
with a property that
$ \lim_{n\rightarrow\infty}P \{L(\alpha_n^0)=
 \inf_{\alpha\in\mathcal{A}_n}L(\alpha) \}=1$; see \citet{Shao1997} for more details. In this case, pursuing consistency is
equivalent to finding the model with the smallest KL loss. However,
$\alpha_n^0$ may not always lead to the smallest KL loss when
$\bolds{\Sigma}_\eta$ has to be estimated, making asymptotic loss
efficiency more difficult to derive. In addition, the possible
inconsistency of $\hat{\bolds{\theta}}(\alpha)$ for
$\alpha\in\mathcal{A}_n$ under the fixed domain asymptotic framework
further complicates the development of asymptotic theory for GIC.

Let $\lambda_{\min}(\mathbf{Q})$ and
$\lambda_{\max}(\mathbf{Q})$ be the smallest and the largest
eigenvalue of a square matrix $\mathbf{Q}$. We impose the following
regularity conditions for model selection:

\begin{longlist}[(C3)]
\item[(C1)]  $\lambda_{\min}(\bolds\Sigma(\bolds\theta))>0$ for all $n$ and $\bolds\theta\in\Theta$, and
\[
\limsup_{n\rightarrow\infty}\sup_{\bolds\theta\in\Theta}\lambda_{\max}
\bigl(\bolds\Sigma^{-1/2}(\bolds\theta)\bolds\Sigma_0\bolds
\Sigma^{-1/2}(\bolds\theta)\bigr)<\infty.
\]

\item[(C2)] For $\alpha\in\mathcal{A}_n\setminus\mathcal{A}_n^0$, there exists
$\bolds\theta_\alpha\in\Theta$, not depending on $n$, such that
\begin{eqnarray*}
\sup_{\alpha\in\mathcal{A}_n\setminus\mathcal{A}_n^0}\biggl |\frac{\ell(\alpha;\hat{\bolds\theta}(\alpha))
-\ell(\alpha;\bolds\theta_\alpha)}{R(\alpha;\bolds\theta_\alpha)
-L_0(\bolds{\theta}_0)} \biggr|&= & o_p(1),
\\
\sup_{\alpha\in\mathcal{A}_n\setminus\mathcal{A}_n^0}\biggl |\frac{L(\alpha;\hat{\bolds\theta}(\alpha))
-L(\alpha;\bolds\theta_\alpha)}{R(\alpha;\bolds\theta_\alpha)
-L_0(\bolds{\theta}_0)} \biggr|& =& o_p(1).
\end{eqnarray*}

Moreover,
\begin{eqnarray*}
\sup_{\alpha\in\mathcal{A}_n^0}\bigl |\ell\bigl(\alpha;\hat{\bolds\theta}(\alpha)\bigr)
-\ell(\alpha;\bolds\theta_0) \bigr| &=& O_p(1),
\\
\sup_{\alpha\in\mathcal{A}_n^0} \bigl|L\bigl(\alpha;\hat{\bolds\theta}(\alpha)\bigr) -L(
\alpha;\bolds\theta_0) \bigr| &=& O_p(1).
\end{eqnarray*}

\item[(C3)] For $\bolds\theta_\alpha$ defined in (C2),
\[
\lim_{n\rightarrow\infty}\sum_{\alpha\in\mathcal{A}_n\setminus\mathcal{A}_n^0}
\frac{1}{(R(\alpha;\bolds\theta_\alpha)-L_0(\bolds{\theta}_0))^q}=0,
\]
for some $q>0$.

\item[(C4)] For $\bolds\theta_\alpha$ defined in (C2),
\[
\lim_{n\rightarrow\infty}\sup_{\alpha\in\mathcal{A}_n \setminus
 \mathcal{A}^{0}_n} \biggl|\frac{\operatorname{tr}
 (\bolds\Sigma_0(\bolds\Sigma^{-1}_{0}-
 \bolds\Sigma^{-1}(\bolds\theta_\alpha))
 \mathbf{M}(\alpha;\bolds\theta_\alpha))}{R(\alpha;\bolds\theta_\alpha)-
 L_0(\bolds{\theta}_0)} \biggr|
=0.
\]
\item[(C5)] For $\bolds\theta_\alpha$ defined in (C2),
\[
\sup_{\alpha\in\mathcal{A}_n \setminus \mathcal{A}^{0}_n} \biggl|\frac{\operatorname{tr} (((\bolds\eta+\bolds\epsilon)(\bolds\eta+\bolds\epsilon)'-\bolds\Sigma_0)
    (\bolds\Sigma^{-1}(\bolds\theta_\alpha)-\bolds\Sigma^{-1}(\bolds\theta_0)) )}    {
R(\alpha;\bolds\theta_\alpha)-L_0(\bolds{\theta}_0)} \biggr|=o_{p}(1).
\]
\end{longlist}

While $L_0(\bolds{\theta}_0)=0$ for a correct spatial covariance model,
we still keep $L_0(\bolds{\theta}_0)$ in (C2)--(C5) because $L_0(\bolds{\theta}_0) \neq 0$
under covariance mis-specification, which will be discussed
in Section~\ref{section:covariance model selection}.
In the rest of this section, we shall assume $\bolds\Sigma_0=\bolds\Sigma(\bolds\theta_0)$,
yielding $L_0(\bolds{\theta}_0)=0$.
Condition (C1), imposing some constraints on the family of covariance
matrices parameterized by $\bolds\theta \in \Theta$, is usually satisfied when $\Theta$ is compact.
Condition (C2) generally holds when $\hat{\bolds\theta}(\alpha)$
converges in probability to some $\bolds\theta_\alpha\in\Theta$,
not necessarily equal to $\bolds\theta_{0}$.
Surprisingly, it can hold even if
$\hat{\bolds\theta}(\alpha)$ does not converge in probability;
see Section~\ref{section:examples} for some examples in which the domain $D$ is fixed with $n$.
Condition (C3) is easily met when
$ |\mathcal{A}_n\setminus\mathcal{A}_n^0 |$ (i.e., the number
of models in $\mathcal{A}_n\setminus\mathcal{A}_n^0$) is bounded and
\[
\min_{\alpha\in\mathcal{A}_n\setminus\mathcal{A}_n^0}\bigl \|\bolds\Sigma^{-1/2}(\bolds
\theta_\alpha)\mathbf{A}(\alpha;\bolds\theta_\alpha) \bolds
\mu_0 \bigr\|^2\rightarrow\infty,
\]
as $n\rightarrow\infty$. Moreover, (C5) can be verified using some moment bounds for quadratic forms in $\bolds\eta+\bolds\epsilon$,
and (C4) is ensured by (C3) when $p_n$ is bounded.

Conditions (C1)--(C5) appear to be natural generalizations
of the conditions used to establish the asymptotic
loss efficiency in usual linear regression models. To see this, note that if $\bolds\Sigma_{0}=\bolds\Sigma(\bolds\theta_{0})$ is known
(or, equivalently, $\Theta=\{\bolds\theta_{0}\}$), then
(C1), (C2), (C4) and (C5) become redundant, and only (C3) is needed, which
corresponds to (A.3) of \citet{Li1987} or (2.6) of \citet{Shao1997}. This is the only assumption needed to
derive the asymptotic loss efficiency
of AIC under model (\ref{setup})
with $\eta(\cdot)=0$, $v^{2}$ known, $\mathbf{s}_{i}=i$; $i=1,\ldots, n$,
and $ |\mathcal{A}_n^0 |\leq 1$.
For more details, see Theorem~1 of \citet{Shao1997}.
On the other hand, when $\bolds\theta_{0}$ is unknown,
(C1), (C2), (C4) and (C5) seem indispensable for dealing with the inherent difficulties in
model selection under (\ref{setup}). That is, the ML estimate of $\bolds\theta$ may not only vary
across candidate models, but may also converge to wrong parameter vectors or have no probability limits.
In the following theorem, these four conditions will be used in conjunction with (C3)
to establish the consistency and
the asymptotic loss efficiency of AIC,
extending Theorem~1 of \citet{Shao1997} to the geostatistical model described in (\ref{setup}) and (\ref{Sigma}).

%th3.1 #&#
\begin{thmm}
Consider the data generated from (\ref{geo data}) and the model
given by (\ref{setup}) and (\ref{Sigma}) with $\bolds{\theta}_0$ being
the true covariance parameter vector [i.e.,
$\operatorname{var}(\mathbf{Z})= \bolds\Sigma(\bolds\theta_0)$].
Suppose that conditions \textup{(C1)--(C5)} are satisfied:

\begin{longlist}[(ii)]
\item[(i)] If $ |\mathcal{A}_n^0 |\leq 1$,
        then $\hat\alpha_2$ is asymptotically loss efficient.
        If, in addition, $ |\mathcal{A}_n^0 |=1$ and $ \limsup_{n\rightarrow\infty}p(\alpha_n^0)<\infty$,
        then $\hat\alpha_2$ is consistent.
\item[(ii)] If $ |\mathcal{A}_n^0 |\geq 2$ for sufficiently large $n$ and either of the following is satisfied
    for some $m>0$,
%
%e3.6 #&#
%e3.7 #&#
\begin{eqnarray}
\lim_{n\rightarrow\infty}\sum_{\alpha\in\mathcal{A}_n^0}
\frac{1}{p^m(\alpha)} &=& 0,\label{infty correct models}
\\
\lim_{n\rightarrow\infty}\sum_{\alpha\in\mathcal{A}_n^0\setminus\{\alpha_n^0\}}
\frac{1}{(p(\alpha)-p(\alpha_n^0))^m} &=& 0. \label{infty correct models 2}
\end{eqnarray}
Then $\hat\alpha_2$ is asymptotically loss efficient. If, in addition,
(\ref{infty correct models 2}) holds and $ \limsup_{n\rightarrow\infty}p(\alpha_n^0)<\infty$,
then $\hat{\alpha}_2$ is consistent.
\end{longlist}
\label{theorem:unknown AIC}
\end{thmm}

\begin{pf}
We begin by showing that
%
%e3.8 #&#
\begin{equation}
\Gamma_2(\alpha) = \nu + 2L(\alpha;\bolds\theta_\alpha) +
o_p\bigl(L(\alpha;\bolds\theta_\alpha)\bigr), \label{loss eff among incorrect models}
\end{equation}
uniformly for
$\alpha\in\mathcal{A}_n\setminus\mathcal{A}_n^0$, where
$\nu=n\log(2\pi)+\log\det(\bolds\Sigma(\bolds\theta_0))+(\bolds\eta+\bolds\epsilon)'\bolds\Sigma^{-1}(\bolds\theta_0)(\bolds\eta+\bolds\epsilon)$
is independent of $\alpha$. By (\ref{unknown GIC}) and (C2), we
have
\begin{eqnarray*}
\Gamma_2(\alpha) &=& -2\ell(\alpha;\bolds\theta_\alpha) + 2p(
\alpha) +o_p\bigl(R(\alpha;\bolds\theta_\alpha)\bigr)
\\
&=& n\log(2\pi) + \log\det\bigl(\bolds\Sigma(\bolds\theta_\alpha)\bigr) +
\mathbf{Z}'\mathbf{A}(\alpha;\bolds\theta_\alpha)'
\bolds\Sigma^{-1}(\bolds\theta_\alpha) \mathbf{A}(\alpha;\bolds
\theta_\alpha)\mathbf{Z}
\\
&&{} +2p(\alpha)+o_p\bigl(R(\alpha;\bolds\theta_\alpha)\bigr)
\\
&=& n\log(2\pi) + \log\det\bigl(\bolds\Sigma(\bolds\theta_\alpha)\bigr) +
\bolds\mu_0'\bolds\Sigma^{-1}(\bolds
\theta_\alpha) \mathbf{A}(\alpha;\bolds\theta_\alpha)\bolds
\mu_0
\\
&&{} +2\bolds\mu_0'\bolds\Sigma^{-1}(\bolds
\theta_\alpha) \mathbf{A}(\alpha;\bolds\theta_\alpha) (\bolds\eta+
\bolds\epsilon)
\\
&&{} +(\bolds\eta+\bolds\epsilon)'\bolds\Sigma^{-1}(\bolds
\theta_\alpha) \mathbf{A}(\alpha;\bolds\theta_\alpha) (\bolds\eta+
\bolds\epsilon) +2p(\alpha)+o_p\bigl(R(\alpha;\bolds
\theta_\alpha)\bigr),
\end{eqnarray*}
uniformly for $\alpha\in\mathcal{A}_n\setminus\mathcal{A}_n^0$.
It follows from (\ref{fn:KL loss}) that
%
%e3.9 #&#
\begin{eqnarray}
\label{loss eff asmp 0} \Gamma_2(\alpha) &=& n\log(2\pi) + \log\det\bigl(\bolds
\Sigma(\bolds\theta_0)\bigr) + (\bolds\eta+\bolds\epsilon)'
\bolds\Sigma^{-1}(\bolds\theta_\alpha) (\bolds\eta+\bolds
\epsilon)
\nonumber\\
&&{} - \operatorname{tr}\bigl(\bolds\Sigma(\bolds\theta_0)\bolds
\Sigma^{-1}(\bolds\theta_\alpha)\bigr) +n + 2L(\alpha;\bolds
\theta_\alpha)
\nonumber\\
&&{} -2(\bolds\eta+\bolds\epsilon)'\bolds\Sigma^{-1}(\bolds
\theta_\alpha) \mathbf{M}(\alpha;\bolds\theta_\alpha) (\bolds\eta+
\bolds\epsilon)
\nonumber\\
&&{} +2\bolds\mu_0'\bolds\Sigma^{-1}(\bolds
\theta_\alpha) \mathbf{A}(\alpha;\bolds\theta_\alpha) (\bolds\eta+
\bolds\epsilon) +2p(\alpha)+o_p\bigl(R(\alpha;\bolds
\theta_\alpha)\bigr)
\nonumber
\\[-8pt]
\\[-8pt]
\nonumber
&=& n\log(2\pi) + \log\det\bigl(\bolds\Sigma(\bolds\theta_0)\bigr) +
(\bolds\eta+\bolds\epsilon)'\bolds\Sigma^{-1}(\bolds
\theta_0) (\bolds\eta+\bolds\epsilon)
\\
&&{} +\operatorname{tr} \bigl(\bigl((\bolds\eta+\bolds\epsilon) (\bolds\eta+\bolds
\epsilon)'-\bolds\Sigma(\bolds\theta_0)\bigr) \bigl(
\bolds\Sigma^{-1}(\bolds\theta_\alpha)-\bolds
\Sigma^{-1}(\bolds\theta_0)\bigr) \bigr)\nonumber
\\
&&{} + 2L(\alpha;\bolds\theta_\alpha) -2(\bolds\eta+\bolds
\epsilon)'\bolds\Sigma^{-1}(\bolds\theta_\alpha)
\mathbf{M}(\alpha;\bolds\theta_\alpha) (\bolds\eta+\bolds\epsilon)+2p(
\alpha)\nonumber
\\
&&{} +2\bolds\mu_0'\bolds\Sigma^{-1}(\bolds
\theta_\alpha) \mathbf{A}(\alpha;\bolds\theta_\alpha) (\bolds\eta+
\bolds\epsilon)+ o_p\bigl(R(\alpha;\bolds\theta_\alpha)\bigr),\nonumber
\end{eqnarray}
uniformly for $\alpha\in\mathcal{A}_n\setminus\mathcal{A}_n^0$.
Therefore, by (C5), for (\ref{loss eff among incorrect
models}) to hold, it suffices to show that
%
%e3.10 #&#
%e3.11 #&#
\begin{eqnarray}
(\bolds\eta+\bolds\epsilon)'\bolds\Sigma^{-1}(\bolds
\theta_\alpha) \mathbf{M}(\alpha;\bolds\theta_\alpha) (\bolds\eta+
\bolds\epsilon)-p(\alpha) &=& o_p\bigl(R(\alpha;\bolds
\theta_\alpha)\bigr), \label{loss eff asmp 2}
\\
\bolds\mu_0'\bolds\Sigma^{-1}(\bolds
\theta_\alpha) \mathbf{A}(\alpha;\bolds\theta_\alpha) (\bolds\eta+
\bolds\epsilon) &=& o_p\bigl(R(\alpha;\bolds\theta_\alpha)
\bigr), \label{loss eff asmp 3}
\end{eqnarray}
uniformly for $\alpha\in\mathcal{A}_n\setminus\mathcal{A}_n^0$,
and
%
%e3.12 #&#
\begin{equation}
\sup_{\alpha\in\mathcal{A}_n\setminus\mathcal{A}_n^0}
\biggl|\frac{L(\alpha;\bolds\theta_\alpha)}{R(\alpha;\bolds\theta_\alpha)}-1\biggr |=o_p(1).
\label{loss eff asmp 1}
\end{equation}

First, we prove (\ref{loss eff asmp 2}). By (C4),
we have
\[
\mathrm{E}\bigl\{(\bolds\eta+\bolds\epsilon)'\bolds
\Sigma^{-1}(\bolds\theta_\alpha) \mathbf{M}(\alpha;\bolds
\theta_\alpha) (\bolds\eta+\bolds\epsilon)\bigr\} - p(\alpha) = o\bigl(R(
\alpha;\bolds\theta_\alpha)\bigr),
\]
uniformly for
$\alpha\in\mathcal{A}_n\setminus\mathcal{A}_n^0$. Let
$c(\alpha)=\operatorname{tr}(\bolds\Sigma(\bolds\theta_0)\bolds\Sigma^{-1}(\bolds\theta_\alpha)\mathbf{M}(\alpha;\bolds\theta_\alpha))/
p(\alpha)$. Then by (\ref{fn:M}) and (C1), $
\limsup_{n\rightarrow\infty}\sup_{\alpha\in\mathcal{A}_n\setminus\mathcal{A}_n^0}c(\alpha)<\infty$.
Thus for (\ref{loss eff asmp 2}) to hold, it suffices to show that
\begin{eqnarray*}
&&(\bolds\eta+\bolds\epsilon)'\bolds\Sigma^{-1}(\bolds
\theta_\alpha) \mathbf{M}(\alpha;\bolds\theta_\alpha) (\bolds\eta+
\bolds\epsilon)-c(\alpha)p(\alpha)\\
&&\qquad = o_p\bigl(R(\alpha;\bolds
\theta_\alpha)\bigr),
\end{eqnarray*}
uniformly for $\alpha\in\mathcal{A}_n\setminus\mathcal{A}_n^0$.
Applying Chebyshev's inequality, we have for any $\varepsilon>0$,
\begin{eqnarray*}
&&P \biggl\{\sup_{\alpha\in\mathcal{A}_n\setminus\mathcal{A}_n^0} \biggl|\frac{(\bolds\eta+\bolds\epsilon)'\bolds\Sigma^{-1}(\bolds\theta_\alpha)
\mathbf{M}(\alpha;\bolds\theta_\alpha)(\bolds\eta+\bolds\epsilon)-c(\alpha)p(\alpha)}{
R(\alpha;\bolds\theta_\alpha)} \biggr| >\varepsilon
\biggr\}
\\
&&\qquad\leq  \sum_{\alpha\in\mathcal{A}_n\setminus\mathcal{A}_n^0} \frac{\mathrm{E} |(\bolds\eta+\bolds\epsilon)'\bolds\Sigma^{-1}(\bolds\theta_\alpha)
\mathbf{M}(\alpha;\bolds\theta_\alpha)(\bolds\eta+\bolds\epsilon)-c(\alpha)p(\alpha) |^{2q}}{
\varepsilon^{2q}R^{2q}(\alpha;\bolds\theta_\alpha)}
\\
&&\qquad\leq  \sum_{\alpha\in\mathcal{A}_n\setminus\mathcal{A}_n^0} \frac{c_1\{\operatorname{tr}(\bolds\Sigma(\bolds\theta_0)\bolds\Sigma^{-1}(\bolds\theta_\alpha)
\mathbf{M}(\alpha;\bolds\theta_\alpha)\bolds\Sigma(\bolds\theta_0)\bolds\Sigma^{-1}(\bolds\theta_\alpha)
\mathbf{M}(\alpha;\bolds\theta_\alpha))\}^q}{
\varepsilon^{2q}R^{2q}(\alpha;\bolds\theta_\alpha)}
\\
&&\qquad\leq \sum_{\alpha\in\mathcal{A}_n\setminus\mathcal{A}_n^0} \frac{c_2 p^q(\alpha)}    {
\varepsilon^{2q}R^{2q}(\alpha;\bolds\theta_\alpha)} \\
&&\qquad\leq \sum
_{\alpha\in\mathcal{A}_n\setminus\mathcal{A}_n^0} \frac{c_3}{\varepsilon^{2q}R^q(\alpha;\bolds\theta_\alpha)},
\end{eqnarray*}
for some constants $c_1,c_2,c_3>0$, where the second
inequality follows from Theorem~2 of \citet{Whittle1960} that
$\mathrm{E}(|\mathbf{y}'\mathbf{A}\mathbf{y}-\mathrm{E}(\mathbf{y}'\mathbf{A}\mathbf{y})|)^{2q}\leq
c_1(\operatorname{tr}(\mathbf{A}^2))^q$ for
$\mathbf{y}=\bolds\Sigma^{-1/2}(\bolds\theta_0)(\bolds\eta+\bolds\epsilon)\sim
N(\mathbf{0},\mathbf{I})$ and
$\mathbf{A}=\bolds\Sigma^{1/2}(\bolds\theta_0)\bolds\Sigma^{-1}(\bolds\theta_\alpha)\mathbf{M}(\alpha;\bolds\theta_\alpha)\bolds\Sigma^{1/2}(\bolds\theta_0)$,
the third inequality follows from (C1), and the last inequality
follows from (C4). Therefore by (C3), we obtain (\ref{loss eff asmp
2}).

Next, we prove (\ref{loss eff asmp 3}). Similar to the proof of (\ref{loss eff asmp 2}),
we have
\begin{eqnarray*}
&& P \biggl\{\sup_{\alpha\in\mathcal{A}_n\setminus\mathcal{A}_n^0} \biggl| \frac{\bolds\mu_0'\bolds\Sigma^{-1}(\bolds\theta_\alpha)
    \mathbf{A}(\alpha;\bolds\theta_\alpha)(\bolds\eta+\bolds\epsilon)}    {
R(\alpha;\bolds\theta_\alpha)} \biggr|>\varepsilon
\biggr\}
\\
&&\qquad\leq  \sum_{\alpha\in\mathcal{A}_n\setminus\mathcal{A}_n^0} \frac{\mathrm{E} |\bolds\mu_0'\bolds\Sigma^{-1}(\bolds\theta_\alpha)
    \mathbf{A}(\alpha;\bolds\theta_\alpha)(\bolds\eta+\bolds\epsilon) |^{2q}}    {
\varepsilon^{2q}R^{2q}(\alpha;\bolds\theta_\alpha)}
\\
&&\qquad\leq  \sum_{\alpha\in\mathcal{A}_n\setminus\mathcal{A}_n^0} \frac{c_4
    (\bolds\mu_0'\bolds\Sigma^{-1}(\bolds\theta_\alpha)
    \mathbf{A}(\alpha;\bolds\theta_\alpha)\bolds\Sigma(\bolds\theta_0)\mathbf{A}(\alpha;\bolds\theta_\alpha)'\bolds\Sigma^{-1}(\bolds\theta_\alpha)\bolds\mu_0)^q}    {
\varepsilon^{2q}R^{2q}(\alpha;\bolds\theta_\alpha)}
\\
&&\qquad\leq  \sum_{\alpha\in\mathcal{A}_n\setminus\mathcal{A}_n^0} \frac{c_5(\bolds\mu_0'\bolds\Sigma^{-1}(\bolds\theta_\alpha)
    \mathbf{A}(\alpha;\bolds\theta_\alpha)\bolds\mu_0)^q}
    {\varepsilon^{2q}R^{2q}(\alpha;\bolds\theta_\alpha)} \\
    &&\qquad\leq \sum
_{\alpha\in\mathcal{A}_n\setminus\mathcal{A}_n^0} \frac{c_6}{\varepsilon^{2q}R^q(\alpha;\bolds\theta_\alpha)},
\end{eqnarray*}
for some constant $c_4,c_5,c_6>0$, where the second
inequality follows from Theorem~2 of \citet{Whittle1960} that
$\mathrm{E}(|\mathbf{a}'\mathbf{y}|)^{2q}\leq c_4(\mathbf{a}'\mathbf{a})^q$ for
$\mathbf{y}=\bolds\Sigma^{-1/2}(\bolds\theta_0)(\bolds\eta+\bolds\epsilon)\sim
N(\mathbf{0},\mathbf{I})$ and
$\mathbf{a}=\bolds\Sigma^{1/2}(\bolds\theta_0)\mathbf{A}(\alpha;\bolds\theta_\alpha)\bolds\Sigma^{-1}(\bolds\theta_\alpha)\bolds{\mu}_0$,
the third inequality follows from~(C1), and the last inequality
follows from (\ref{fn:KL risk}). Therefore by (C3), we obtain~(\ref{loss eff asmp 3}).

It remains to prove (\ref{loss eff asmp 1}). By (\ref{fn:KL
loss}) and (\ref{fn:KL risk}), for $\alpha\in\mathcal{A}_n\setminus\mathcal{A}_n^0$,
\begin{eqnarray*}
L(\alpha;\bolds\theta_\alpha) - R(\alpha;\bolds\theta_\alpha) &=&
(\bolds\eta+\bolds\epsilon)'\bolds\Sigma^{-1}(\bolds
\theta_\alpha) \mathbf{M}(\alpha;\bolds\theta_\alpha) (\bolds\eta+
\bolds\epsilon)
\\
&&{} -\operatorname{tr}\bigl(\bolds\Sigma^{-1}(\bolds
\theta_\alpha) \mathbf{M}(\alpha;\bolds\theta_\alpha)\bolds\Sigma(
\bolds\theta_0)\bigr).
\end{eqnarray*}
It follows from (C1), (C3) and an
argument similar to one used to prove (\ref{loss eff asmp 2}) that
\begin{eqnarray*}
&&\sup_{\alpha\in\mathcal{A}_n\setminus\mathcal{A}_n^0} \biggl|
\frac{L(\alpha;\bolds\theta_\alpha)}{R(\alpha;\bolds\theta_\alpha)}-1 \biggr|\\
&&\qquad =\sup_{\alpha\in\mathcal{A}_n\setminus\mathcal{A}_n^0}
\biggl|  \frac{(\bolds\eta+\bolds\epsilon)'\bolds\Sigma^{-1}(\bolds\theta_\alpha)
    \mathbf{M}(\alpha;\bolds\theta_\alpha)(\bolds\eta+\bolds\epsilon)}    {
R(\alpha;\bolds\theta_\alpha)}
 \\
 &&\hspace*{83pt}{}-\frac{\operatorname{tr}(\bolds\Sigma^{-1}(\bolds\theta_\alpha)
    \mathbf{M}(\alpha;\bolds\theta_\alpha)\bolds\Sigma(\bolds\theta_0))
    }{R(\alpha;\bolds\theta_\alpha)} \biggr|=o_p(1).
\end{eqnarray*}
 This gives (\ref{loss eff asmp 1}). Thus (\ref{loss eff
among incorrect models}) is established.

(i) If $ |\mathcal{A}_n^0 |=0$, it follows from
(\ref{loss eff among incorrect models}), (\ref{loss eff asmp 1}) and (C2) that
$\hat{\alpha}_2$ is asymptotically loss efficient. If
$ |\mathcal{A}_n^0 |= 1$ and
$ \lim_{n\rightarrow\infty}p(\alpha_n^0)=\infty$,
by (\ref{loss eff among incorrect models}), to
show the asymptotic loss efficiency of $\hat\alpha_2$, it suffices to show that
%
%e3.13 #&#
\begin{equation}
\Gamma_2(\alpha) = \nu + 2L(\alpha;\bolds\theta_0) +
o_p\bigl(L(\alpha;\bolds\theta_0)\bigr); \qquad \alpha\in
\mathcal{A}_n^0. \label{uniform result 1}
\end{equation}
%  uniformly for $\alpha\in\mathcal{A}_n^0$.
By (\ref{log like fun alpha}), (\ref{fn:KL loss 0}) and (C2),
%
%e3.14 #&#
\begin{eqnarray}
\label{eq:AIC} \Gamma_2(\alpha) &=& -2\ell(\alpha;\bolds
\theta_0) + 2p(\alpha)+ O_p(1)\nonumber
\\
&=& \nu -2\bigl\{(\bolds\eta+\bolds\epsilon)'\bolds
\Sigma^{-1}(\bolds\theta_0) \mathbf{M}(\alpha;\bolds
\theta_0) (\bolds\eta+\bolds\epsilon)-p(\alpha)\bigr\}
\\
&&{} + 2L(\alpha;\bolds\theta_0)+ O_p(1);\qquad \alpha\in
\mathcal{A}_n^0.\nonumber %  \alpha\in\mathcal{A}_n^0.
\end{eqnarray}
Therefore, by (\ref{fn:KL loss 0}), (\ref{fn:KL risk 0}) and
an argument similar to that used to prove (\ref{loss eff among incorrect
models}), we have
%
%e3.15 #&#
%e3.16 #&#
\begin{eqnarray}
\biggl|\frac{
(\bolds\eta+\bolds\epsilon)'\bolds\Sigma^{-1}(\bolds\theta_0)
\mathbf{M}(\alpha_n^0;\bolds\theta_0)(\bolds\eta+\bolds\epsilon)
-p(\alpha_n^0)}{p(\alpha_n^0)} \biggr| &=& o_p(1), \label{uniform result 1-1}
\\
\biggl|\frac{L(\alpha_n^0;\bolds\theta_0)}{R(\alpha_n^0;\bolds\theta_0)}-1 \biggr| &=& o_p(1). \label{uniform result 1-3}
\end{eqnarray}
These together with (\ref{eq:AIC}) give (\ref{uniform result 1}).
If $ |\mathcal{A}_n^0 |=1$ and $ \limsup_{n\rightarrow\infty}p(\alpha_n^0)<\infty$,
then the consistency and the asymptotical loss efficiency are ensured by
%
%e3.17 #&#
%e3.18 #&#
\begin{eqnarray}
L\bigl(\alpha;\hat{\bolds\theta}(\alpha)\bigr) - L\bigl(\alpha_n^0;
\hat{\bolds\theta}\bigl(\alpha_n^0\bigr)\bigr)
&\mathop{\rightarrow}\limits^{P}& \infty, \label{AIC:eq2}
\\
\Gamma_2(\alpha)-\Gamma_2\bigl(\alpha_n^0
\bigr) &\mathop{\rightarrow}\limits^{P}& \infty, \label{AIC:eq3}
\end{eqnarray}
uniformly for $\alpha\in\mathcal{A}_n\setminus\{\alpha_n^0\}$, as
$n\rightarrow\infty$. First, (\ref{AIC:eq2}) follows from
%
%e3.19 #&#
%e3.20 #&#
%e3.21 #&#
\begin{eqnarray}
\label{smallest loss} L\bigl(\alpha_n^0;\hat{\bolds\theta}\bigl(
\alpha_n^0\bigr)\bigr)& =& L\bigl(\alpha_n^0;
\bolds\theta_0\bigr) + O_p(1)\nonumber
\\
&=& \tfrac{1}{2}(\bolds\eta+\bolds\epsilon)'\bolds
\Sigma^{-1}(\bolds\theta_0) \mathbf{M}\bigl(
\alpha_n^0;\bolds\theta_0\bigr) (\bolds\eta+
\bolds\epsilon) +O_p(1)
\\
&=& o_p\bigl(L\bigl(\alpha;\hat{\bolds\theta}(\alpha)\bigr)\bigr),\nonumber
\end{eqnarray}
uniformly for $\alpha\in\mathcal{A}_n\setminus\{\alpha_n^0\}$, where the first equality follows from (C2), the second
equality follows from (\ref{fn:KL loss}) and the last equality
follows from (\ref{loss eff asmp 1}), (C2), (C3) and $  \limsup_{n\rightarrow\infty}p(\alpha_n^0)<\infty$.
It remains to prove (\ref{AIC:eq3}). By (\ref{eq:AIC}), we have
%
%e3.22 #&#
\begin{eqnarray}
\label{AIC for correct model0} \Gamma_2\bigl(\alpha_n^0\bigr)
&=& \nu -(\bolds\eta+\bolds\epsilon)'\bolds\Sigma^{-1}(
\bolds\theta_0) \mathbf{M}\bigl(\alpha_n^0;
\bolds\theta_0\bigr) (\bolds\eta+\bolds\epsilon)+2p\bigl(
\alpha_n^0\bigr)+ O_p(1)
\nonumber
\\[-8pt]
\\[-8pt]
\nonumber
&=& \nu + o_p\bigl(L\bigl(\alpha;\hat{\bolds\theta}(\alpha)\bigr)
\bigr),
\end{eqnarray}
uniformly for $\alpha\in\mathcal{A}_n\setminus\{\alpha_n^0\}$, where the last
equality follows from an argument similar to that used to prove (\ref{smallest loss}).
This together with (\ref{loss eff among incorrect models}) implies
(\ref{AIC:eq3}). This completes the proof of (i).

(ii) First, suppose that (\ref{infty correct models}) is satisfied.
In view of (\ref{loss eff among incorrect models}), it suffices to show that~(\ref{uniform result 1}) holds
uniformly for $\alpha\in\mathcal{A}_n^0$.
Similarly to the proofs of (\ref{uniform result 1-1}) and (\ref{uniform result 1-3}), we only need to show that
%
%e3.23 #&#
%e3.24 #&#
\begin{eqnarray}
\sup_{\alpha\in\mathcal{A}_n^0} \biggl|\frac{
(\bolds\eta+\bolds\epsilon)'\bolds\Sigma^{-1}(\bolds\theta_0)
\mathbf{M}(\alpha;\bolds\theta_0)(\bolds\eta+\bolds\epsilon)
-p(\alpha)}{p(\alpha)} \biggr| &=& o_p(1),
\label{uniform result 1-2}
\\
\sup_{\alpha\in\mathcal{A}_n^0}
 \biggl|\frac{L(\alpha;\bolds\theta_0)}{R(\alpha;\bolds\theta_0)}-1 \biggr| &=& o_p(1).
\label{uniform result 1-4}
\end{eqnarray}

By an argument similar to that used to prove (\ref{loss eff asmp
2}), we have
\begin{eqnarray*}
&& P \biggl\{ \sup_{\alpha\in\mathcal{A}_n^0} \biggl|\frac{
(\bolds\eta+\bolds\epsilon)'\bolds\Sigma^{-1}(\bolds\theta_0)
\mathbf{M}(\alpha;\bolds\theta_0)(\bolds\eta+\bolds\epsilon)-p(\alpha)}
{p(\alpha)}\biggr |>\varepsilon
\biggr\}
\\
 &&\qquad\leq \sum_{\alpha\in\mathcal{A}_n^0}\frac{\mathrm{E} |(\bolds\eta+\bolds\epsilon)'\bolds\Sigma^{-1}(\bolds\theta_0)
\mathbf{M}(\alpha;\bolds\theta_0)(\bolds\eta+\bolds\epsilon)-p(\alpha) |^{2m}}{\varepsilon^{2m}p^{2m}(\alpha)}
\\
&&\qquad \leq \sum
_{\alpha\in\mathcal{A}_n^0}\frac{c_7}{\varepsilon^{2m}p^m(\alpha)},
\end{eqnarray*}
for some constant $c_7>0$, as $n\rightarrow\infty$. This together with
 (\ref{fn:KL loss 0}), (\ref{fn:KL risk 0}) and (\ref{infty correct models})
 gives~(\ref{uniform result 1-2}) and (\ref{uniform result 1-4}).
Therefore, (\ref{uniform result 1}) holds uniformly for $\alpha\in\mathcal{A}_n^0$.

Finally, suppose that (\ref{infty correct models 2}) is satisfied. If
$ \lim_{n\rightarrow\infty}p(\alpha_n^0)=\infty$, it
implies (\ref{infty correct models}) and hence $\hat\alpha_2$ is
asymptotically loss efficient. If
$ \limsup_{n\rightarrow\infty}p(\alpha_n^0)<\infty$, by
(\ref{AIC:eq2}) and (\ref{AIC:eq3}), it remains to show that
%
%e3.25 #&#
%e3.26 #&#
\begin{eqnarray}
L\bigl(\alpha;\hat{\bolds\theta}(\alpha)\bigr)- L\bigl(\alpha_n^0;
\hat{\bolds\theta}\bigl(\alpha_n^0\bigr)\bigr)
&\mathop{\rightarrow}\limits^{P}& \infty, \label{AIC:eq4}
\\
\Gamma_2(\alpha)-\Gamma_2\bigl(\alpha_n^0
\bigr) &\mathop{\rightarrow}\limits^{P}& \infty, \label{AIC:eq5}
\end{eqnarray}
uniformly for
$\alpha\in\mathcal{A}_n^0\setminus\{\alpha_n^0\}$, as
$n\rightarrow\infty$. First, we prove (\ref{AIC:eq4}). By
(\ref{fn:KL risk 0}) and (C2),
%
%e3.27 #&#
\begin{eqnarray}
\label{AIC:eq6} &&L\bigl(\alpha;\hat{\bolds\theta}(\alpha)\bigr)-L\bigl(
\alpha_n^0;\hat{\bolds\theta}\bigl(\alpha_n^0
\bigr)\bigr)\nonumber
\\
&&\qquad= L(\alpha;\bolds\theta_0)-L\bigl(\alpha_n^0;
\bolds\theta_0\bigr) + O_p(1)
\nonumber
\\[-8pt]
\\[-8pt]
\nonumber
&&\qquad= \tfrac{1}{2}(\bolds\eta+\bolds\epsilon)'\bolds
\Sigma^{-1}(\bolds\theta_0) \bigl\{\mathbf{M}(\alpha;\bolds
\theta_0)-\mathbf{M}\bigl(\alpha_n^0;\bolds
\theta_0\bigr)\bigr\}(\bolds\eta+\bolds\epsilon)+O_p(1)
\\
&&\qquad= \tfrac{1}{2} \bigl(p(\alpha)-p\bigl(\alpha_n^0
\bigr) \bigr)+ o_p\bigl(p(\alpha)-p\bigl(\alpha_n^0
\bigr)\bigr),\nonumber
\end{eqnarray}
uniformly for
$\alpha\in\mathcal{A}_n^0\setminus\{\alpha_n^0\}$, where the last
equality follows from
\begin{eqnarray*}
&&\sup_{\alpha\in\mathcal{A}_n^0\setminus\{\alpha_n^0\}}\biggl|\frac{
(\bolds\eta+\bolds\epsilon)'\bolds\Sigma^{-1}(\bolds\theta_0)
\mathbf{M}(\alpha;\bolds\theta_0)(\bolds\eta+\bolds\epsilon)-p(\alpha)}{p(\alpha)-p(\alpha_n^0)}
\\
&&\hspace*{21pt}\qquad{} -\frac{
(\bolds\eta+\bolds\epsilon)'\bolds\Sigma^{-1}(\bolds\theta_0)
\mathbf{M}(\alpha_n^0;\bolds\theta_0)(\bolds\eta+\bolds\epsilon)
-p(\alpha_n^0)}{p(\alpha)-p(\alpha_n^0)}\biggr | =o_p(1),
\end{eqnarray*}
which can be obtained in a way similar to the proof of
(\ref{uniform result 1-1}). This together with~(\ref{infty correct
models 2}) gives (\ref{AIC:eq4}). Next, we prove (\ref{AIC:eq5}). By
(\ref{eq:AIC}) and (\ref{AIC:eq6}), we have
\begin{eqnarray*}
\Gamma_2(\alpha)-\Gamma_2\bigl(\alpha_n^0
\bigr) &=& 2L\bigl(\alpha;\hat{\bolds\theta}(\alpha)\bigr)-2L\bigl(
\alpha_n^0;\hat{\bolds\theta}\bigl(\alpha_n^0
\bigr)\bigr) +o_p\bigl(p(\alpha)-p\bigl(\alpha_n^0
\bigr)\bigr)
\\
&=& p(\alpha)-p\bigl(\alpha_n^0\bigr)+o_p
\bigl(p(\alpha)-p\bigl(\alpha_n^0\bigr)\bigr),
\end{eqnarray*}
 uniformly for $\alpha\in\mathcal{A}_n^0\setminus\{\alpha_n^0\}$.
This together with (\ref{infty correct models 2}) gives
(\ref{AIC:eq5}). This completes the proof of (ii).
\end{pf}

%re3.1 #&#
\begin{rem}
When $\bolds\theta=\bolds\theta_0$ is known, Theorem~\ref{theorem:unknown
AIC} reduces to the standard asymptotic theory of AIC in linear
regression; see Theorem~1 of \citet{Shao1997}. In this case, (C1), (C2), (C4)
and (C5) are not needed.
\end{rem}

%re3.2 #&#
\begin{rem}
Although Theorem~\ref{theorem:unknown AIC} only obtains
the consistency of $\hat\alpha_2$ under $ \limsup_{n\rightarrow\infty}p(\alpha_n^0)<\infty$,
the consistency result can be extended to
$ \lim_{n\rightarrow\infty}p(\alpha_n^0)=\infty$
if
$p(\alpha_n^0)=o ( \inf_{\alpha\in\mathcal{A}_n\setminus\mathcal{A}_n^0}R(\alpha;\bolds\theta_0) )$.
\end{rem}

%re3.3 #&#
\begin{rem}
When $ |\mathcal{A}_n^0 |\geq 2$, AIC is generally not able to identify
$\alpha_n^0$ almost surely. A heavier penalty ${\tau_n}$ of GIC (e.g., BIC) is
needed for consistency.
\end{rem}

%th3.2 #&#
\begin{thmm}
Consider the data generated from (\ref{geo data}) and the model
given by (\ref{setup}) and (\ref{Sigma}) with $\bolds{\theta}_0$ being
the true covariance parameter vector [i.e.,
$\operatorname{var}(\mathbf{Z})= \bolds\Sigma(\bolds\theta_0)$]. Suppose
that \textup{(C1)--(C5)} are satisfied. In addition, suppose that
$ \lim_{n\rightarrow\infty}{\tau_n}=\infty$, and for
$\bolds\theta_\alpha$ defined in \textup{(C2)},
%
%e3.28 #&#
\begin{equation}
\lim_{n\rightarrow\infty}\sup_{\alpha\in\mathcal{A}_n\setminus\mathcal{A}_n^0}\frac{\tau_n p_n}{R(\alpha;\bolds\theta_\alpha)} =0.
\label{eq:C6}
\end{equation}

\begin{longlist}[(ii)]
\item[(i)] If $ |\mathcal{A}_n^0 |=0$, then $\hat\alpha_{\tau_n}$ is asymptotically loss efficient.
\item[(ii)] If $ |\mathcal{A}_n^0 |\geq 1$ and
%
%e3.29 #&#
\begin{equation}
\lim_{n\rightarrow\infty}\sum_{\alpha\in\mathcal{A}_n^0}
\frac{1}{p^m(\alpha)}<\infty, \label{cond for GIC}
\end{equation}
for some $m>0$, then $\hat\alpha_{\tau_n}$ is
consistent.
\end{longlist}
\label{theorem:unknown GIC}
\end{thmm}

\begin{pf}
(i) By (\ref{loss eff among incorrect models}) and
(\ref{eq:C6}), we have
%
%e3.30 #&#
\begin{equation}
\Gamma_{\tau_n}(\alpha) =\nu + 2L(\alpha;\bolds\theta_\alpha) +
o_p\bigl(L(\alpha;\bolds\theta_\alpha)\bigr), \label{loss eff among incorrect models 2}
\end{equation}
uniformly for $\alpha\in\mathcal{A}_n\setminus\mathcal{A}_n^0$. Thus by (\ref{loss eff asmp 1}) and (C2),
$\hat\alpha_{\tau_n}$ is asymptotically loss efficient.

(ii) By (\ref{log like fun alpha}) and (C2), we have for
$\alpha\in\mathcal{A}_n^0$,
%
%e3.31 #&#
\begin{eqnarray}
\label{GIC:eq2} \Gamma_{\tau_n}(\alpha)& =& -2\ell(\alpha;\bolds
\theta_0) + {\tau_n} p(\alpha)+ O_p(1)
\nonumber
\\[-8pt]
\\[-8pt]
\nonumber
&=& \nu - (\bolds\eta+\bolds\epsilon)'\bolds\Sigma^{-1}(
\bolds\theta_0) \mathbf{M}(\alpha;\bolds\theta_0) (\bolds
\eta+\bolds\epsilon) +{\tau_n} p(\alpha) + O_p(1),
\end{eqnarray}
where $\nu$ is defined in (\ref{loss eff among incorrect models}).
By (\ref{cond for GIC}) and an argument similar to that used to prove~(\ref{uniform result 1-2}), we have
%
%e3.32 #&#
\begin{equation}
\sup_{\alpha\in\mathcal{A}_n^0} \biggl| \frac{(\bolds\eta+\bolds\epsilon)'
\bolds\Sigma^{-1}(\bolds\theta_0)\mathbf{M}(\alpha;\bolds\theta_0)
(\bolds\eta+\bolds\epsilon)-p(\alpha)}{{\tau_n} p(\alpha)} \biggr|=o_p(1).
\label{GIC:eq1}
\end{equation}
This and (\ref{GIC:eq2}) give
%
%e3.33 #&#
\begin{equation}
\Gamma_{\tau_n}(\alpha) = \nu + (\tau_n-1) p(\alpha) +
o_p\bigl({\tau_n} p(\alpha)\bigr), \label{GIC:eq3}
\end{equation}
uniformly for $\alpha\in\mathcal{A}_n^0$. Thus
%
%e3.34 #&#
\begin{equation}
\lim_{n\rightarrow\infty}P\bigl\{\hat{\alpha}_{\tau_n}\in
\mathcal{A}_n^0\setminus\bigl\{\alpha_n^0
\bigr\}\bigr\} =0. \label{GIC:eq4}
\end{equation}
By (\ref{eq:C6}), (\ref{loss eff among incorrect
models 2}) and (\ref{GIC:eq3}), we have
\[
\min_{\alpha\in\mathcal{A}_n\setminus\mathcal{A}_n^0} \Gamma_{\tau_n}(\alpha) -
\Gamma_{\tau_n}\bigl(\alpha_n^0\bigr)
\mathop{\rightarrow}\limits^{P}\infty,
\]
as $n\rightarrow\infty$. This together with
(\ref{GIC:eq4}) implies that $\hat\alpha_{\tau_n}$ is consistent. This
completes the proof.
\end{pf}

Unlike the KL loss function in usual linear regression
models, $L(\alpha, \hat{\bolds\theta}(\alpha))$ does not necessarily
have the minimum at $\alpha=\alpha^{0}_{n}$, and hence selection
consistency may not lead to asymptotic loss efficiency in
geostatistical regression models. Nevertheless, when
$\bolds\theta=\bolds\theta_0$ is known, Theorem~\ref{theorem:unknown GIC}
reduces to the standard asymptotic theory of GIC in linear
regression [see Theorem~2 of \citet{Shao1997}], in which selection
consistency is known to imply asymptotic loss efficiency. This
property continues to hold if $\hat{\bolds{\theta}}(\alpha)$ in
(\ref{unknown GIC}), and (\ref{loss efficiency}) is replaced by a common estimate
$\hat{\bolds{\theta}}$, independent of $\alpha$. Then for
$\alpha\in\mathcal{A}_n^0\setminus\{\alpha_n^0\}$,
\[
L(\alpha;\hat{\bolds\theta})-L\bigl(\alpha_n^0;\hat{
\bolds\theta}\bigr) =(\bolds\eta+\bolds\epsilon)'\bolds
\Sigma^{-1}(\hat{\bolds\theta}) \bigl(\mathbf{M}(\alpha;\hat{\bolds
\theta})-\mathbf{M}\bigl(\alpha_n^0;\hat{\bolds\theta}
\bigr)\bigr) (\bolds\eta+\bolds\epsilon)\geq 0,
\]
 almost surely.

%co3.1 #&#
\begin{coro}\label{coro:unknown GIC}
Consider the data generated from (\ref{geo data}) and the model
defined in (\ref{setup}) and (\ref{Sigma}) with $\bolds{\theta}_0$
being the true covariance parameter vector [i.e.,
$\operatorname{var}(\mathbf{Z})= \bolds\Sigma(\bolds\theta_0)$]. Suppose
that \textup{(C1)--(C5)} are satisfied with $\hat{\bolds{\theta}}(\alpha)$ and
$\bolds{\theta}_\alpha$ in \textup{(C2)--(C5)} being replaced by
$\hat{\bolds\theta}$ and a constant vector $\bolds{\theta}_c\in\Theta$,
independent of $\alpha$. Let $\hat{\alpha}_{\tau_n}$ be the model
selected by a modified GIC criterion with
$\hat{\bolds{\theta}}(\alpha)$ in (\ref{unknown GIC}) being replaced by
$\hat{\bolds\theta}$. In addition, suppose that
$ \lim_{n\rightarrow\infty}{\tau_n}=\infty$, and
$ \lim_{n\rightarrow\infty} \sup_{\alpha\in\mathcal{A}_n\setminus\mathcal{A}_n^0}\frac{\tau_n
p_n} {R(\alpha;\bolds\theta_c)}=0$.

\begin{longlist}[(ii)]
\item[(i)] If $ |\mathcal{A}_n^0 |=0$, then $\hat\alpha_{\tau_n}$ is asymptotically loss efficient in the sense that
    $L(\hat{\alpha}_{\tau_n};\hat{\bolds\theta})
     / \inf_{\alpha\in\mathcal{A}_n}L(\alpha;\hat{\bolds\theta})
    \mathop{\rightarrow}\limits^{P}1$, as $n\rightarrow\infty$.
\item[(ii)] If $ |\mathcal{A}_n^0 |\geq 1$ and (\ref{cond for GIC}) holds, then
        $\hat\alpha_{\tau_n}$ is consistent and asymptotically
    loss efficient in the sense that
    $L(\hat{\alpha}_{\tau_n};\hat{\bolds\theta})
     / \inf_{\alpha\in\mathcal{A}_n}L(\alpha;\hat{\bolds\theta})
    \mathop{\rightarrow}\limits^{P}1$, as $n\rightarrow\infty$.
\end{longlist}
\end{coro}

%s4 #&#
\section{Variable selection under an incorrect covariance model}
\label{section:covariance model selection}

In this section, we establish the asymptotic theory of GIC for
variable selection, when the covariance model is mis-specified with
$\bolds{\Sigma}_0\neq\bolds{\Sigma}(\bolds{\theta}_0)$, yielding $L_0(\bolds{\theta}_0)\neq 0$.
To ensure that the asymptotic optimality of GIC for $\bolds{\Sigma}_0=\bolds{\Sigma}(\bolds{\theta}_0)$
carries over to this case, we need a stronger condition in place of (C4):

\begin{enumerate}
\item[(C4$'$)] For $\bolds\theta_\alpha$ defined in (C2),
\[
\lim_{n\rightarrow\infty}\sup_{\alpha\in\mathcal{A}_n\setminus\mathcal{A}_n^0} \frac{p_n}{R(\alpha;\bolds\theta_\alpha)-L_0(\bolds{\theta}_0)} =0.
\]
\end{enumerate}

%th4.1 #&#
\begin{thmm}\label{theorem:unknown AIC 2}
Consider the data generated from (\ref{geo data}) and the model
given by (\ref{setup}) and (\ref{Sigma}). Suppose that the
conditions \textup{(C1)--(C3), (C4$'$)} and \textup{(C5)} are satisfied:

\begin{longlist}[(ii)]
\item[(i)] If $ |\mathcal{A}_n^0 |\leq 1$, then
$\hat\alpha_2$ is asymptotically loss efficient. If
$ |\mathcal{A}_n^0 |= 1$, then $\hat\alpha_2$ is
consistent.
\item[(ii)] If $ |\mathcal{A}_n^0 |\geq 2$ for sufficient large $n$,
$ |\mathcal{A}_n^0 |^q=o(L_0(\bolds\theta_0))$ for some $q>0$,
and
%
%e4.1 #&#
\begin{equation}
\lim_{n\rightarrow\infty}\frac{p_n}{L_0(\bolds\theta_0)} =0, \label{eq:C7}
\end{equation}
then $\hat\alpha_2$ is asymptotically loss efficient.
\end{longlist}
\end{thmm}

\begin{pf}
Let $L^*(\alpha;\bolds\theta_\alpha)=L(\alpha;\bolds\theta_\alpha)-L_0(\bolds{\theta}_0)$;
$\alpha\in\mathcal{A}_n\setminus\mathcal{A}^0_n$. We begin by showing that
%
%e4.2 #&#
\begin{equation}
\Gamma_2(\alpha) = \nu+ 2L^*(\alpha;\bolds\theta_\alpha) +
o_p\bigl(L^*(\alpha;\bolds\theta_\alpha)\bigr),
\label{loss eff among incorrect models 3}
\end{equation}
uniformly for
$\alpha\in\mathcal{A}_n\setminus\mathcal{A}_n^0$, where
$\nu$ is defined in (\ref{loss eff among incorrect models}) and is
independent of $\alpha$. By an argument similar to that used to prove
(\ref{loss eff asmp 0}), we have
\begin{eqnarray*}
\Gamma_2(\alpha) &=& n\log(2\pi) + \log\det(\bolds
\Sigma_0) + n -\operatorname{tr}\bigl(\bolds\Sigma_0\bolds
\Sigma^{-1}(\bolds\theta)\bigr)
\\
&&{} + (\bolds\eta+\bolds\epsilon)'\bolds\Sigma^{-1}(\bolds
\theta_0) (\bolds\eta+\bolds\epsilon)
\\
&&{} +\operatorname{tr} \bigl(\bigl((\bolds\eta+\bolds\epsilon) (\bolds\eta+\bolds
\epsilon)'-\bolds\Sigma_0\bigr) \bigl(\bolds
\Sigma^{-1}(\bolds\theta_\alpha)-\bolds\Sigma^{-1}(
\bolds\theta_0)\bigr) \bigr)
\\
&&{} +2L(\alpha;\bolds\theta_\alpha) -2(\bolds\eta+\bolds
\epsilon)'\bolds\Sigma^{-1}(\bolds\theta_\alpha)
\mathbf{M}(\alpha;\bolds\theta_\alpha) (\bolds\eta+\bolds\epsilon)+2p(
\alpha)
\\
&&{} +2\bolds\mu_0'\bolds\Sigma^{-1}(\bolds
\theta_\alpha) \mathbf{A}(\alpha;\bolds\theta_\alpha) (\bolds\eta+
\bolds\epsilon)+ o_p\bigl(R(\alpha;\bolds\theta_\alpha)\bigr)
\\
&=& \nu +\operatorname{tr} \bigl(\bigl((\bolds\eta+\bolds\epsilon) (\bolds\eta+
\bolds\epsilon)'-\bolds\Sigma_0\bigr) \bigl(\bolds
\Sigma^{-1}(\bolds\theta_\alpha)-\bolds\Sigma^{-1}(
\bolds\theta_0)\bigr) \bigr)
\\
&&{} + 2L^*(\alpha;\bolds\theta_\alpha)-2(\bolds\eta+\bolds
\epsilon)'\bolds\Sigma^{-1}(\bolds\theta_\alpha)
\mathbf{M}(\alpha;\bolds\theta_\alpha) (\bolds\eta+\bolds\epsilon)+2p(
\alpha)
\\
&&{} +2\bolds\mu_0'\bolds\Sigma^{-1}(\bolds
\theta_\alpha) \mathbf{A}(\alpha;\bolds\theta_\alpha) (\bolds\eta+
\bolds\epsilon)+ o_p\bigl(R(\alpha;\bolds\theta_\alpha)\bigr),
\end{eqnarray*}
uniformly for
$\alpha\in\mathcal{A}_n\setminus\mathcal{A}_n^0$. Hence by (C5) and
an argument similar to that used to prove (\ref{loss eff among incorrect
models}), for (\ref{loss eff among incorrect models 3}) to hold, it
suffices to show that
\begin{eqnarray*}
(\bolds\eta+\bolds\epsilon)'\bolds\Sigma^{-1}(\bolds
\theta_\alpha) \mathbf{M}(\alpha;\bolds\theta_\alpha) (\bolds\eta+
\bolds\epsilon)-p(\alpha) &=& o_p\bigl(R(\alpha;\bolds
\theta_\alpha)-L_0(\bolds{\theta}_0)\bigr),
\\
\bolds\mu_0'\bolds\Sigma^{-1}(\bolds
\theta_\alpha) \mathbf{A}(\alpha;\bolds\theta_\alpha) (\bolds\eta+
\bolds\epsilon) &=& o_p\bigl(R(\alpha;\bolds\theta_\alpha)-L_0(
\bolds{\theta}_0)\bigr),
\end{eqnarray*}
uniformly for
$\alpha\in\mathcal{A}_n\setminus\mathcal{A}_n^0$, and
%
%e4.3 #&#
\begin{equation}
\sup_{\alpha\in\mathcal{A}_n
\setminus\mathcal{A}_n^0} \biggl|\frac{L^*(\alpha;\bolds\theta_\alpha)}
{R(\alpha;\bolds\theta_\alpha)-L_0(\bolds{\theta}_0)}-1 \biggr|=o_p(1).
\label{loss eff asmp 6}
\end{equation}
The above three equations follow from arguments similar to those used
 to prove~(\ref{loss eff asmp 2})--(\ref{loss eff asmp 1}).

(i) Clearly, (\ref{loss eff among incorrect models 3}) implies
(\ref{loss eff among incorrect models}). Therefore, if
$ |\mathcal{A}_n^0 |= 0$, it follows from
(\ref{loss eff asmp 6}) and (C2) that $\hat{\alpha}_2$ is
asymptotically loss efficient. On the other hand, if
\mbox{$ |\mathcal{A}_n^0 |=1$}, it suffices to show
(\ref{AIC:eq2}) and (\ref{AIC:eq3}). First, we prove
(\ref{AIC:eq2}). By (C3), (C4$'$) and an argument similar to that used to prove
(\ref{smallest loss}), we have
$L^*(\alpha_n^0;\hat{\bolds\theta}(\alpha_n^0))=o_p(L^*(\alpha;\hat{\bolds\theta}(\alpha)))$,
uniformly for $\alpha\in\mathcal{A}_n\setminus\{\alpha_n^0\}$. Next,
we prove (\ref{AIC:eq3}). By (C3), (C4$'$) and an argument similar to
that used to prove (\ref{AIC for correct model0}), we have
$\Gamma_2(\alpha_n^0)=\nu^* +
o_p(L^*(\alpha;\hat{\bolds\theta}(\alpha)))$, uniformly for
$\alpha\in\mathcal{A}_n\setminus\{\alpha_n^0\}$. This together with
(\ref{loss eff among incorrect models 3}) implies (\ref{AIC:eq3}), and
hence the proof of (i) is complete.

(ii) In view of (\ref{loss eff among incorrect models}), it suffices to show
that
%
%e4.4 #&#
\begin{equation}
\Gamma_2(\alpha) = \nu^* + 2L(\alpha;\bolds\theta_0) +
o_p\bigl(L(\alpha;\bolds\theta_0)\bigr), \label{AIC for correct model-3}
\end{equation}
uniformly for $\alpha\in\mathcal{A}_n^0$, where
$\nu^*=\nu-2L_0(\bolds\theta_0)$ with $\nu$ being defined in (\ref{loss eff
among incorrect models}). By an argument similar to that used to prove (\ref{eq:AIC}),
we have
%
%e4.5 #&#
\begin{eqnarray}
\label{loss eff among incorrect models 7} \Gamma_2(\alpha) &=& \nu^* -2\bigl\{(\bolds\eta+\bolds
\epsilon)'\bolds\Sigma^{-1}(\bolds\theta_0)
\mathbf{M}(\alpha;\bolds\theta_0) (\bolds\eta+\bolds\epsilon)-p(
\alpha)\bigr\}
\nonumber
\\[-8pt]
\\[-8pt]
\nonumber
&&{} + 2L(\alpha;\bolds\theta_0)+ O_p(1); \qquad\alpha\in
\mathcal{A}_n^0.
\end{eqnarray}
Therefore, by an argument similar to that used to prove (\ref{uniform
result 1}), we only need to show that
%
%e4.6 #&#
\begin{equation}
(\bolds\eta+\bolds\epsilon)'\bolds\Sigma^{-1}(\bolds
\theta_0) \mathbf{M}(\alpha;\bolds\theta_0) (\bolds\eta+
\bolds\epsilon)-p(\alpha) = o_p\bigl(L(\alpha;\bolds
\theta_0)\bigr), \label{loss eff asmp 7}
\end{equation}
uniformly for $\alpha\in\mathcal{A}_n^0$ and
%
%e4.7 #&#
\begin{equation}
\sup_{\alpha\in\mathcal{A}_n^0}
\biggl|\frac{L(\alpha;\bolds{\theta}_0)}{R(\alpha;\bolds{\theta}_0)}-1 \biggr| =o_p(1).
\label{loss eff asmp 8}
\end{equation}

First, we prove (\ref{loss eff asmp 7}). Clearly, by
(\ref{fn:M}) and (C1), we have
%
%e4.8 #&#
\begin{equation}
\mathrm{E}\bigl((\bolds\eta+\bolds\epsilon)'\bolds
\Sigma^{-1}(\bolds\theta_0) \mathbf{M}(\alpha;\bolds
\theta_0) (\bolds\eta+\bolds\epsilon)\bigr)=c(\alpha)p(\alpha),
\label{AIC:eq1}
\end{equation}
where
$ \limsup_{n\rightarrow\infty}\sup_{\alpha\in\mathcal{A}_n}c(\alpha)<\infty$.
Hence by (\ref{fn:KL risk}) and (\ref{eq:C7}),
$c(\alpha)p(\alpha)-p(\alpha)=o(R(\alpha;\bolds\theta_0))$ uniformly
for $\alpha\in\mathcal{A}_n^0$. It remains to show that
\[
(\bolds\eta+\bolds\epsilon)'\bolds\Sigma^{-1}(\bolds
\theta_0) \mathbf{M}(\alpha;\bolds\theta_0) (\bolds\eta+
\bolds\epsilon) - c(\alpha)p(\alpha) = o_p\bigl(R(\alpha;\bolds
\theta_0)\bigr),
\]
uniformly for $\alpha\in\mathcal{A}_n^0$. Applying
Chebyshev's inequality, we have for any $\varepsilon>0$,
\begin{eqnarray*}
&& P \biggl\{\sup_{\alpha\in\mathcal{A}_n^0}\biggl |\frac{(\bolds\eta+\bolds\epsilon)'\bolds\Sigma^{-1}(\bolds\theta_0)
\mathbf{M}(\alpha;\bolds\theta_0)(\bolds\eta+\bolds\epsilon)-c(\alpha)p(\alpha)}{
R(\alpha;\bolds\theta_0)} \biggr| >\varepsilon
\biggr\}
\\
&&\qquad\leq  \sum_{\alpha\in\mathcal{A}_n^0} \frac{\mathrm{E} |(\bolds\eta+\bolds\epsilon)'\bolds\Sigma^{-1}(\bolds\theta_0)
\mathbf{M}(\alpha;\bolds\theta_0)(\bolds\eta+\bolds\epsilon)-c(\alpha)p(\alpha) |^{2m}}{
\varepsilon^{2m}R^{2m}(\alpha;\bolds\theta_0)}
\\
&&\qquad\leq  \sum_{\alpha\in\mathcal{A}_n^0} \frac{c_1\{\operatorname{tr}(\bolds\Sigma_0\bolds\Sigma^{-1}(\bolds\theta_0)
\mathbf{M}(\alpha;\bolds\theta_0)\bolds\Sigma_0\bolds\Sigma^{-1}(\bolds\theta_0)
\mathbf{M}(\alpha;\bolds\theta_0))\}^m}{
\varepsilon^{2m}R^{2m}(\alpha;\bolds\theta_0)}
\\
&&\qquad\leq \sum_{\alpha\in\mathcal{A}_n^0} \frac{c_2 p^m(\alpha)}    {
\varepsilon^{2m}L_0^{2m}(\bolds\theta_0)} \leq \sum
_{\alpha\in\mathcal{A}_n^0} \frac{c_3}{\varepsilon^{2m}L_0^m(\bolds\theta_0)},
\end{eqnarray*}
where the second-to-last equality follows from (C1) and
$R(\alpha;\bolds\theta_0)\geq L_0(\bolds\theta_0)$, for
$\alpha\in\mathcal{A}_n$, and the last equality follows from
(\ref{eq:C7}). Taking $m=1/q$, we obtain (\ref{loss eff asmp 7}). Next, we prove (\ref{loss eff asmp 8}). By (\ref{fn:KL
loss}), (\ref{fn:KL risk}) and (\ref{AIC:eq1}), we have for
$\alpha\in\mathcal{A}_n^0$,
\[
L(\alpha;\bolds\theta_0)-R(\alpha;\bolds\theta_0) =
\tfrac{1}{2} \bigl\{(\bolds\eta+\bolds\epsilon)'\bolds
\Sigma^{-1}(\bolds\theta_0) \mathbf{M}(\alpha;\bolds
\theta_0) (\bolds\eta+\bolds\epsilon) - c(\alpha)p(\alpha) \bigr\},
\]
where
$ \limsup_{n\rightarrow\infty}\sup_{\alpha\in\mathcal{A}_n^0}c(\alpha)<\infty$.
Thus (\ref{loss eff asmp 8}) follows from an argument similar to that used to prove
(\ref{loss eff asmp 7}).\vspace*{1pt} Thus we obtain (\ref{AIC for
correct model-3}). This completes the proof.
\end{pf}

%th4.2 #&#
\begin{thmm}\label{theorem:unknown GIC 2}
Under the setup of Theorem~\ref{theorem:unknown AIC 2}, suppose that\break
$ \lim_{n\rightarrow\infty}\tau_n =\infty$, and
%
%e4.9 #&#
\begin{equation}
\lim_{n\rightarrow\infty}\sup_{\alpha\in\mathcal{A}_n\setminus\mathcal{A}_n^0} \frac{\tau_n p_n}{R(\alpha;\bolds\theta_\alpha)-L_0(\bolds{\theta}_0)} = 0.
\label{eq:C8}
\end{equation}

\begin{longlist}[(ii)]
\item[(i)] If $ |\mathcal{A}_n^0 |=0$, then $\hat\alpha_{\tau_n}$ is asymptotically loss efficient.
\item[(ii)] If $ |\mathcal{A}_n^0 |\geq 1$, $ |\mathcal{A}_n^0 |^q=o(L_0(\bolds\theta_0))$ for
some $q>0$, and (\ref{cond for GIC}) is satisfied, then
$\hat\alpha_{\tau_n}$ is consistent and asymptotically loss
efficient.
\end{longlist}
\end{thmm}

\begin{pf}
(i) By (\ref{loss eff among incorrect models 3}) and
(\ref{eq:C8}), we have
% \begin{eqnarray*}
$\Gamma_{\tau_n}(\alpha)
= \nu + 2L^*(\alpha;\bolds\theta_\alpha) +\break o_p(L^*(\alpha; \bolds\theta_\alpha))$,
% \label{loss eff among incorrect models 4}
% \end{eqnarray*}
uniformly for
$\alpha\in\mathcal{A}_n\setminus\mathcal{A}_n^0$, and hence
%
%e4.10 #&#
\begin{equation}
\Gamma_{\tau_n}(\alpha) = \nu^* + 2L(\alpha;\bolds\theta_\alpha)
+ o_p\bigl(L(\alpha;\bolds\theta_\alpha)\bigr),
\label{loss eff among incorrect models 5}
\end{equation}
uniformly for
$\alpha\in\mathcal{A}_n\setminus\mathcal{A}_n^0$. In addition, (\ref{loss eff asmp 6}) gives
%
%e4.11 #&#
\begin{equation}
\sup_{\alpha\in\mathcal{A}_n\setminus\mathcal{A}_n^0}
\biggl |\frac{L(\alpha;\bolds\theta_\alpha)}{R(\alpha;\bolds\theta_\alpha)}-1 \biggr|=o_p(1).
\label{loss eff among incorrect models 6}
\end{equation}
These together with (C2) imply that $\hat\alpha_{\tau_n}$ is asymptotically loss efficient.

(ii) First, we prove the asymptotic loss efficiency of
$\hat\alpha_{\tau_n}$. By (\ref{loss eff asmp 8}) and
(\ref{loss eff among incorrect models 6}), we have
%
%e4.12 #&#
\begin{equation}
\sup_{\alpha\in\mathcal{A}_n}
 \biggl|\frac{L(\alpha;\bolds{\theta}_0)}{R(\alpha;\bolds{\theta}_0)}-1 \biggr| =o_p(1).
\label{loss eff asmp 9}
\end{equation}
By (\ref{eq:C8}) and an argument similar to that used to prove
(\ref{AIC for correct model-3}), we have
\[
\Gamma_{\tau_n}(\alpha) = \nu^* + 2 L(\alpha;\bolds\theta_0)
+ o_p\bigl(L(\alpha;\bolds\theta_0)\bigr),
\]
 uniformly for $\alpha\in\mathcal{A}_n^0$. This together
with (\ref{loss eff among incorrect models 5}), (\ref{loss
eff asmp 9}) and (C2) implies that $\hat\alpha_{\tau_n}$ is
asymptotically loss efficient.

Next, we prove the consistency of $\hat\alpha_{\tau_n}$. By
(\ref{loss eff among incorrect models 7}) and (\ref{AIC:eq1}), we
have for $\alpha\in\mathcal{A}_n^0$,
%
%e4.13 #&#
\begin{eqnarray}
\label{GIC2:eq2} \Gamma_{\tau_n}(\alpha) &=& \nu - \bigl\{(\bolds\eta+\bolds
\epsilon)'\bolds\Sigma^{-1}(\bolds\theta_0)
\mathbf{M}(\alpha;\bolds\theta_0) (\bolds\eta+\bolds\epsilon)-c(
\alpha)p(\alpha)\bigr\}
\nonumber
\\[-8pt]
\\[-8pt]
\nonumber
&&{} +\bigl(\tau_n-c(\alpha)\bigr) p(\alpha) + o_p\bigl(
\tau_n p(\alpha)\bigr).
\end{eqnarray}
By (\ref{cond for GIC}) and an argument similar to that used to prove (\ref{GIC:eq1}), we have
\[
\sup_{\alpha\in\mathcal{A}_n^0} \biggl| \frac{(\bolds\eta+\bolds\epsilon)'
\bolds\Sigma^{-1}(\bolds\theta_0)\mathbf{M}(\alpha;\bolds\theta_0)
(\bolds\eta+\bolds\epsilon)-c(\alpha)p(\alpha)}{{\tau_n} p(\alpha)} \biggr|=o_p(1).
\]
Hence by (\ref{GIC2:eq2}),
%
%e4.14 #&#
\begin{equation}
\Gamma_{\tau_n}(\alpha) = \nu + \bigl(\tau_n-c(\alpha)\bigr)
p(\alpha) + o_p\bigl({\tau_n} p(\alpha)\bigr),
\label{GIC2:eq3}
\end{equation}
uniformly for $\alpha\in\mathcal{A}_n^0$. Thus we obtain
(\ref{GIC:eq4}). In addition, by (\ref{eq:C8}), (\ref{loss eff among
incorrect models 5}) and~(\ref{GIC2:eq3}),
\[
\min_{\alpha\in\mathcal{A}_n\setminus\mathcal{A}_n^0} \Gamma_{\tau_n}(\alpha) -
\Gamma_{\tau_n}\bigl(\alpha_n^0\bigr)
\mathop{\rightarrow}\limits^{P}\infty,
\]
as $n\rightarrow\infty$. This together with
(\ref{GIC:eq4}) implies that $\hat\alpha_{\tau_n}$ is consistent.
This completes the proof.
\end{pf}

%re4.1 #&#
\begin{rem}
Recall that in (ii) of Theorem~\ref{theorem:unknown GIC}, asymptotic
loss efficiency of GIC is generally not satisfied, unless
$\hat{\bolds{\theta}}(\alpha)$'s are replaced by a common estimate. In
contrast, in (ii) of Theorem~\ref{theorem:unknown GIC 2}, we have,
from (\ref{fn:KL loss}) and an argument similar to that used to prove
(\ref{loss eff asmp 7}) that $L(\alpha;\bolds\theta_0)=L_0(\bolds\theta_0)+
o_p(L_0(\bolds\theta_0))$, uniformly for $\alpha\in\mathcal{A}_n^0$,
which leads to
\[
\frac{L(\alpha;\hat{\bolds\theta}(\alpha))}{ \min_{\alpha'\in\mathcal{A}_n}
L(\alpha';\hat{\bolds\theta}(\alpha'))}  \mathop{\rightarrow}\limits^{P}1,
\]
for any $\alpha\in\mathcal{A}_n^0$, indicating that the
asymptotic loss efficiency can be achieved for any correct model.
\end{rem}

%----------------------------------------------------

%s5 #&#
\section{Examples}
\label{section:examples}

In this section, we provide some specific examples for GIC that
satisfy regularity conditions (C1)--(C5). Throughout this
section, we assume that $p_n=p$, $\mathcal{A}_n=\mathcal{A}$,
$\mathcal{A}_n^0=\mathcal{A}^0$ and $\alpha_n^0=\alpha^0$ are fixed, and give
proofs of the theoretical results in the supplemental
material [Chang, Huang and Ing (\citeyear{supp})].

%s5.1 #&#
\subsection{One-dimensional examples}
\label{one dim exp model}

First, we consider spatial models in the one-dimensional space with
$D=[0,n^\delta]\subseteq\mathbb{R}$; $\delta\in[0,1)$. We assume the
exponential covariance model for $\eta(\cdot)$,
%
%e5.1 #&#
\begin{equation}
\operatorname{cov}\bigl(\eta(s),\eta\bigl(s^*\bigr)\bigr) = \sigma^2
\exp\bigl(-\kappa\bigl|s-s^*\bigr|\bigr);\qquad s, s^*\in D, \label{exp cov fun 1 dim}
\end{equation}
where $\sigma^2>0$ is the variance parameter, and
$\kappa>0$ is a spatial dependence parameter. We also assume that
the data are uniformly sampled at $s_i=i n^{-(1-\delta)}$;
$i=1,\ldots,n, s_i\in D$. Clearly, $\delta=0$ corresponds to the
fixed domain asymptotic framework with $D=[0,1]$, and  a larger
$\delta$ corresponds to a faster growth rate of the domain. Note
that $\sigma^2\kappa$ is often referred to as a microergodic
parameter under fixed domain asymptotics [\citet{Stein1999}].

The following proposition allows us to replace
(C1)--(C5) in Theorems \ref{theorem:unknown AIC}
and~\ref{theorem:unknown GIC} by simpler conditions.

%pr5.1 #&#
\begin{pro}\label{pro:1 dim}
Consider $\bolds\Sigma(\bolds\theta)$ in (\ref{Sigma}), where
$\bolds\Sigma_\eta$ is given by (\ref{exp cov fun 1 dim}) and
$s_i=in^{-(1-\delta)}$; $i=1,\ldots,n$, for some $\delta\in[0,1)$.
Let $\bolds\theta=(v^2,\sigma^2,\kappa)'$. Then for any compact set
$\Theta\subseteq(0,\infty)^3$ and any
$\bolds\theta_0=(v_0^2,\sigma_0^2,\kappa_0)'\in\Theta$,
%
%e5.2 #&#
\begin{eqnarray}
\label{proposition:bound eigen} 0 &< & \liminf_{n\rightarrow\infty} \inf_{\bolds\theta\in\Theta}
\lambda_{\min} \bigl( \bolds\Sigma^{-1/2}(\bolds\theta)\bolds
\Sigma(\bolds\theta_0) \bolds\Sigma^{-1/2}(\bolds\theta)
\bigr)
\nonumber
\\[-8pt]
\\[-8pt]
\nonumber
&\leq & \limsup_{n\rightarrow\infty}\sup_{\bolds\theta\in\Theta}
\lambda_{\max} \bigl(\bolds\Sigma^{-1/2}(\bolds\theta) \bolds
\Sigma(\bolds\theta_0)\bolds\Sigma^{-1/2}(\bolds\theta)
\bigr)< \infty.
\end{eqnarray}
\end{pro}

\begin{pf}
The proof follows directly from Proposition~2.1 of Chang, Huang and Ing (\citeyear{Chang2013}).
\end{pf}

% We have the following theorem for GIC.

%th5.1 #&#
\begin{thmm}\label{theorem of GIC in exp model}
Consider the data generated from (\ref{geo data}) and the model
given by (\ref{setup}) and (\ref{Sigma}) with $\bolds{\theta}_0$ being
the true covariance parameter vector [i.e.,
$\operatorname{var}(\mathbf{Z})= \bolds\Sigma(\bolds\theta_0)$]. Assume the
setup of Proposition~\ref{pro:1 dim} with $\delta\in(0,1)$. Suppose
that $\hat{\bolds\theta}(\alpha)\mathop{\rightarrow}\limits^{P}\bolds\theta_\alpha$ for
some $\bolds\theta_\alpha\in\Theta$; $\alpha\in\mathcal{A}$, and
%
%e5.3 #&#
\begin{equation}
\min_{\alpha\in\mathcal{A}\setminus\mathcal{A}^0}
R(\alpha;\bolds\theta_\alpha)\rightarrow
\infty, \label{assump of risk}
\end{equation}
as $n\rightarrow\infty$. Then $\hat\alpha_2$ is
asymptotically loss efficient if $ |\mathcal{A}^0 |\leq 1$. In
addition, suppose that
$ \lim_{n\rightarrow\infty}{\tau_n}=\infty$ and
$ {\tau_n}=o (\min_{\alpha\in\mathcal{A}\setminus\mathcal{A}^0}R(\alpha;\bolds\theta_\alpha) )$.

\begin{longlist}[(ii)]
\item[(i)] If $ |\mathcal{A}^0 |=0$, then $\hat\alpha_{\tau_n}$ is asymptotically loss efficient.
\item[(ii)] If $ |\mathcal{A}^0 |\geq 1$, then $\hat\alpha_{\tau_n}$ is consistent.
\end{longlist}
\end{thmm}

%re5.1 #&#
\begin{rem}
The assumption,
$\hat{\bolds\theta}(\alpha)\mathop{\rightarrow}\limits^{P}\bolds\theta_\alpha$;
$\alpha\in\mathcal{A}$, is generally satisfied under the increasing
domain asymptotic framework, and is guaranteed to hold when
$R(\alpha;\bolds\theta_0)=o(n^{\delta})$, for all
$\alpha\in\mathcal{A}\setminus\mathcal{A}^0$; see Theorem~2.3 of
Chang, Huang and Ing (\citeyear{Chang2013}). In fact, as given by Theorems
\ref{theorem white-noise gic unknown}--\ref{theorem poly gic
unknown}, the assumption continues to hold even if
$R(\alpha;\bolds\theta_0)>cn^{\delta}$ for
$\alpha\in\mathcal{A}\setminus\mathcal{A}^0$ and some constant $c>0$.
\end{rem}

Although the theorem is established under the increasing domain
asymptotic framework, the theorem remains valid in some situations
even when $\hat{\bolds\theta}(\alpha)$ fails to converge for some
$\alpha\in\mathcal{A}$ under the fixed domain asymptotic framework
with \mbox{$\delta=0$}. As mentioned at the end of Section~\ref{subsec:geo model}, our asymptotic results of GIC are still valid
for random $\mathbf{X}$. In what follows, we provide three
examples based on different classes of regressors that are
either random or fixed. We derive the consistency of GIC
not only for $\delta\in(0,1)$ but also for $\delta=0$ without
requiring the regularity conditions. The three examples below can be
seen to have increasing degrees of smoothness in space, leading to
different conditions to ensure the consistency of GIC.

%ex5.1 #&#
\begin{exm}[(White-noise processes)]\label{exm:white-noise}
Consider $p$ regressors, $x_j(\cdot)$; $j=1,\ldots,p$,
% sampled at $s_i=in^{-(1-\delta)}$; $i=1,\ldots,n$, where $x_1(\cdot),\ldots,x_p(\cdot)$ are
generated from independent white-noise processes with
\[
x_j(s)\sim N\bigl(0,v_j^2\bigr);\qquad s\in
\bigl[0,n^\delta\bigr], j=1,\ldots,p,
\]
for some $\delta\in[0,1)$, where $v_j^2>0$; $j=1,\ldots,p$.
\end{exm}

%ex5.2 #&#
\begin{exm}[(Spatially dependent processes)]\label{exm:exp var}
Consider $p$ regressors, $x_j(\cdot)$; $j=1,\ldots,p$, generated from
independent zero-mean Gaussian spatial processes with covariance
functions
\[
\operatorname{cov}\bigl(x_j(s),x_j\bigl(s'
\bigr)\bigr) = \sigma_j^2 \exp \bigl(-
\kappa_j\bigl|s-s'\bigr| \bigr);\qquad s,s'\in
\bigl[0,n^\delta\bigr],
\]
for some $\delta\in[0,1)$, where
$\sigma_j^2,\kappa_j>0$; $j=1,\ldots,p$.
\end{exm}

%ex5.3 #&#
\begin{exm}[(Monomials)]\label{exm:poly}
Consider $p$ regressors, $x_j(\cdot)$; $j=1,\ldots,p$,
\[
x_j(s)=n^{-\delta j}s^j; \qquad s\in\bigl[0,n^\delta
\bigr],
\]
for some $\delta\in[0,1)$. Note that a scaling factor
$n^{-\delta j}$ is introduced to standardize $x_j(\cdot)$ so that
$ \frac{1}{n^\delta}\int_0^{n^\delta}(x_j(s)-\bar{x}_j)^2
\,ds$ does not depend on $n$, where
$\bar{x}_j= \frac{1}{n^\delta}\int_0^{n^\delta}x_j(s)\,ds$.
\end{exm}

%th5.2 #&#
\begin{thmm}\label{theorem white-noise gic unknown}
Consider the model defined in (\ref{setup}) with the white-noise
regressors given by Example~\ref{exm:white-noise}. Suppose that
$\mathbf{Z}\sim N(\mathbf{X}\bolds\beta_0,\bolds\Sigma(\bolds{\theta}_0))$
conditional on~$\mathbf{X}$, where
$\bolds\beta_0=(\beta_{0,0},\ldots,\beta_{0,p})'\in\mathbb{R}^{p+1}$
and
$\bolds\theta_0=(v_0^2,\sigma_0^2,\kappa_0)'\in\Theta\subseteq(0,\infty)^3$
are constant vectors, and $\bolds{\Sigma}(\bolds{\theta}_0)$ is given by
Proposition~\ref{pro:1 dim} for some $\delta\in[0,1)$. Assume that
$\Theta$ is compact and
\[
\bolds{\theta}_0+ \biggl( \sum_{j\in\alpha^0\setminus\alpha}
\beta_{0,j}^2v_j^2, 0, 0
\biggr)'\in\Theta; \qquad \alpha\in\mathcal{A}.
\]
If $ \lim_{n\rightarrow\infty}{\tau_n}=\infty$
and ${\tau_n}=o(n)$, then $ \lim_{n\rightarrow\infty}
P \{\hat{\alpha}_{\tau_n} = \alpha^0 \}=1$.
\end{thmm}

%re5.2 #&#
\begin{rem}
Theorem~\ref{theorem white-noise gic unknown} assumes $\mathcal{A}^0\neq\varnothing$.
Suppose that $\mu_0(\cdot)$ has an additional unobserved
term $\zeta(\cdot)$, which is also a white-noise process,
%
%e5.4 #&#
\begin{equation}
\mu_0(s) = \beta_{0,0} + \sum
_{j=1}^p\beta_{0,j}x_j(s) +
\zeta(s);\qquad s\in D, \label{assump of mean}
\end{equation}
and hence $ |\mathcal{A}^0 |=0$. Then by Theorem~\ref{theorem of GIC in exp model} and an argument similar to that in proof
of Theorem~\ref{theorem white-noise gic unknown} for $\delta=0$, GIC
is also asymptotically loss efficient for
$\delta\in[0,1)$, provided that
$ \lim_{n\rightarrow\infty}{\tau_n}=\infty$ and
${\tau_n}=o(n)$.
\end{rem}

%th5.3 #&#
\begin{thmm}\label{theorem exp var gic unknown}
Consider the model defined in (\ref{setup}) with the spatially
dependent regressors given by Example~\ref{exm:exp var}. Suppose
that $\mathbf{Z}\sim N(\mathbf{X}\bolds\beta_0,\bolds\Sigma(\bolds{\theta}_0))$
conditional on $\mathbf{X}$, where
$\bolds\beta_0=(\beta_{0,0},\ldots,\beta_{0,p})'\in\mathbb{R}^{p+1}$
and
$\bolds\theta_0=(v_0^2,\sigma_0^2,\kappa_0)'\in\Theta\subseteq(0,\infty)^3$
are constant vectors, and $\bolds{\Sigma}(\bolds{\theta}_0)$ is given by
Proposition~\ref{pro:1 dim} with $\delta\in[0,1)$. Assume that
$\Theta$ is compact and
$\bolds{\theta}_0+ (0,  \sum_{j\in\alpha^0\setminus\alpha}\beta_{0,j}^2\sigma_j^2,
\kappa_\alpha^* )'\in\Theta$ for any $\alpha\in\mathcal{A}$,
where $  \kappa_\alpha^*
= (\sigma_0^2+\sum_{j\in\alpha^0\setminus\alpha}\beta_{0,j}^2\sigma_j^2 )^{-1}
 (\sum_{j\in\alpha^0\setminus\alpha}\beta_{0,j}^2\sigma_j^2(\kappa_j-\kappa_0) )$.
If $ \lim_{n\rightarrow\infty}{\tau_n}=\infty$ and
${\tau_n}=o(n^{(1+\delta)/2})$, then
$ \lim_{n\rightarrow\infty} P \{\hat{\alpha}_{\tau_n}
= \alpha^0 \}=1$.
\end{thmm}

%re5.3 #&#
\begin{rem}
Theorem~\ref{theorem exp var gic unknown} assumes $\mathcal{A}^0\neq\varnothing$. Suppose that $\mu_0(\cdot)$ is given by~(\ref{assump of mean}), where $\zeta(\cdot)$ is an unobserved
spatial dependent
process given in Example~\ref{exm:exp var}.
Then by Theorem~\ref{theorem of GIC in exp model} and an argument similar
to that in proof of Theorem~\ref{theorem exp var gic unknown} for
$\delta=0$, GIC is also asymptotically loss efficient
for $\delta\in[0,1)$, provided that
$ \lim_{n\rightarrow\infty}{\tau_n}=\infty$ and
${\tau_n}=o(n^{(1+\delta)/2})$.
\end{rem}

%th5.4 #&#
\begin{thmm}\label{theorem
poly gic unknown}
Consider the model defined in (\ref{setup}) with the monomial
regressors given by Example~\ref{exm:poly}. Suppose that $\mathbf{Z}\sim
N(\mathbf{X}\bolds\beta_0,\bolds\Sigma(\bolds{\theta}_0))$, where
$\bolds\beta_0=(\beta_{0,0},\ldots,\beta_{0,p})'\in\mathbb{R}^{p+1}$
and
$\bolds\theta_0=(v_0^2,\sigma_0^2,\kappa_0)'\in\Theta\subseteq(0,\infty)^3$
are constant vectors, and $\bolds{\Sigma}(\bolds{\theta}_0)$ is given by
Proposition~\ref{pro:1 dim} with $\delta\in(0,1)$. Assume that
$\mathcal{A}=\{\varnothing,\{1\}, \{1,2\},\ldots,\{1,\ldots,p\}\}$, $\Theta$ is
compact,\vspace*{2pt} and
$\bolds{\theta}_0+ (0, \gamma(k),  -(\sigma_0^2+\gamma(k))^{-1}\gamma(k)\kappa_0 )'$
$\in\Theta$; $k=0,1,\ldots,p$, where
$\gamma(k)=\bolds\beta_0'\mathbf{V}_{p,p}\bolds\beta_0
-\bolds\beta_0'\mathbf{V}_{p,k}\mathbf{V}_{k,k}^{-1}\*\mathbf{V}_{k,p}
\bolds\beta_0$ and
$\mathbf{V}_{k,p}= ( \frac{1}{i+j-1} )_{(k+1)\times (p+1)}$.
If $ \lim_{n\rightarrow\infty}{\tau_n}=\infty$ and
${\tau_n}=o(n^\delta)$, then $ \lim_{n\rightarrow\infty}
P \{\hat{\alpha}_{\tau_n} = \alpha^0 \}=1$.
\end{thmm}

%re5.4 #&#
\begin{rem}
Theorem~\ref{theorem poly gic unknown} assumes
$\mathcal{A}^0\neq\varnothing$. Suppose that $\mu_0(\cdot)$ is given by~(\ref{assump of mean}), where $\zeta(s)=n^{-\delta k}s^k$; $s\in D$,
is an unobserved function with $k>p$. Then by Theorem~\ref{theorem of GIC in exp model}, GIC can be shown to be
asymptotically loss efficient for $\delta\in(0,1)$, provided that
$ \lim_{n\rightarrow\infty}{\tau_n}=\infty$ and
${\tau_n}=o(n^\delta)$.
\end{rem}

The results of Theorems \ref{theorem white-noise gic
unknown}--\ref{theorem poly gic unknown} show that the consistency of
GIC depends on not only the smoothness of regressors in space but
also the growth rate of the domain. Evidently, GIC is more
difficult to identify the true model when the candidate regressors
are smoother in space. Although there exists ${\tau_n}$ such that
GIC is consistent for either white-noise regressors or spatially
dependent regressors under the fixed domain asymptotic
framework, interestingly, as shown in the next theorem,
consistent polynomial order selection turns out not possible when
the true model has at least one nonzero regression coefficient and
$ |\mathcal{A}^0 |\geq 2$ under the fixed domain asymptotic
framework.

%th5.5 #&#
\begin{thmm}[(Inconsistency)]\label{thmm:inconsistent:poly} Consider the same setup as in
Theorem~\ref{theorem poly gic unknown}, except that $\delta=0$:
\begin{itemize}[(ii)]
\item[(i)] If $ \lim_{n\rightarrow\infty}{\tau_n}=\infty$, then
  $ \lim_{n\rightarrow\infty}
  P\{\hat{\alpha}_{\tau_n} = \{\varnothing\}\}=1$.
\item[(ii)] If $\alpha^0\neq\{\varnothing\}$ and $ \liminf_{n\rightarrow\infty}\tau_n>0$, then
  $ \lim_{n\rightarrow\infty}
  P\{\hat{\alpha}_{\tau_n} = \alpha^0\}<1$.
\end{itemize}
\end{thmm}

%s5.2 #&#
\subsection{A two-dimensional exponential model}
\label{two dim exp model}

Consider the multiplicative exponential covariance model
%
%e5.5 #&#
\begin{equation}
\operatorname{cov}\bigl(\eta(\mathbf{s}),\eta\bigl(\mathbf{s}^*\bigr)\bigr) =
\sigma^2\exp\bigl(-\kappa \bigl\{\bigl|s_1-s_1^*\bigr|+\bigl|s_2-s_2^*\bigr|
\bigr\}\bigr), \label{exp cov fun 2 dim}
\end{equation}
parameterized by $\sigma^2>0$ and $\kappa>0$, where
$\mathbf{s}=(s_1,s_2)$ and $\mathbf{s}^*=(s_1^*,s_2^*)\in
D=[0,n^{\delta/2}]^2\subseteq \mathbb{R}^2; \delta\in[0,1)$.
Clearly, $\delta=0$ corresponds to the fixed domain asymptotic
framework with $D=[0,1]^2$, and  a larger $\delta$ corresponds to a
faster growth rate of the domain.

Similarly to the one-dimensional case, we first prove
(\ref{proposition:bound eigen}), which is the key to show (C1)--(C5).

%pr5.2 #&#
\begin{pro}\label{pro:2 dim}
Consider $\bolds\Sigma(\bolds\theta)$ in (\ref{Sigma}) with
$\bolds\Sigma_\eta$ given by (\ref{exp cov fun 2 dim}), $v^2=0$, and
$\mathbf{s}_k= (i m^{-(1-\delta)},j m^{-(1-\delta)} )$;
$k={i+(j-1)}m$; $i,j=1,\ldots,m$, for some integer $m=n^{1/2}$,
where $\delta\in[0,1)$. Let $\bolds\theta=(\sigma^2,\kappa)'$. Then
(\ref{proposition:bound eigen}) holds for any compact set
$\Theta\subseteq(0,\infty)^2$ and any
$\bolds\theta_0=(\sigma_0^2,\kappa_0)'\in\Theta$.
\end{pro}

\begin{pf}
Write
%
%e5.6 #&#
\begin{equation}
\bolds\Sigma(\bolds\theta) = \sigma^2\mathbf{B}(\bolds\theta)\otimes
\mathbf{B}(\bolds\theta), \label{Sigma 2 dim}
\end{equation}
where $\mathbf{B}(\bolds\theta)= (\rho^{|i-j|} )_{m\times
m}$ and $\rho=\exp(-\kappa m^{-(1-\delta)})$. By
(\ref{Sigma 2 dim}),
\begin{eqnarray*}
&&\lambda_{\max} \bigl(\bolds\Sigma^{-1/2}(\bolds\theta)\bolds
\Sigma(\bolds\theta_0)\bolds\Sigma^{-1/2}(\bolds\theta)\bigr)
\\
&&\qquad\leq  \frac{\sigma_0^2}{\sigma^2}\lambda_{\max} \bigl( \bigl(\mathbf{B}(\bolds
\theta_0)\otimes\mathbf{B}(\bolds\theta_0)\bigr) \bigl(
\mathbf{B}^{-1}(\bolds\theta)\otimes\mathbf{B}^{-1}(\bolds
\theta)\bigr) \bigr)
\\
& &\qquad=\frac{\sigma_0^2}{\sigma^2}\lambda_{\max} \bigl(\bigl(\mathbf{B}(\bolds
\theta_0)\mathbf{B}^{-1}(\bolds\theta)\bigr) \otimes\bigl(
\mathbf{B}(\bolds\theta_0)\mathbf{B}^{-1}(\bolds\theta)
\bigr) \bigr)
\\
&&\qquad= \frac{\sigma_0^2}{\sigma^2}\lambda_{\max}^2 \bigl(\bigl(
\mathbf{B}(\bolds\theta_0)\mathbf{B}^{-1}(\bolds\theta)
\bigr) \bigr) <\infty,
\end{eqnarray*}
where the last inequality follows from Proposition~2.1 of
Chang, Huang and Ing (\citeyear{Chang2013}). This gives the last inequality of
(\ref{proposition:bound eigen}). The proof for the first
inequality of~(\ref{proposition:bound eigen}) is analogous and
omitted. This completes the proof.
\end{pf}

%th5.6 #&#
\begin{thmm}\label{theorem of GIC in exp model 2}
Consider the data generated from (\ref{geo data}), the model given
by~(\ref{setup}) and (\ref{Sigma}) and the setup of Proposition~\ref{pro:2 dim} with $\delta\in[0,1)$. Suppose that
$\hat{\bolds\theta}(\alpha)\mathop{\rightarrow}\limits^{P}\bolds\theta_\alpha$ for some
$\bolds\theta_\alpha\in\Theta$; $\alpha\in\mathcal{A}$, and
(\ref{assump of risk}) holds. Then $\hat\alpha_2$ is asymptotically
loss efficient if $ |\mathcal{A}^0 |\leq 1$. In addition,
suppose that $ \lim_{n\rightarrow\infty}{\tau_n}=\infty$
and
$ {\tau_n}=o (\min_{\alpha\in\mathcal{A}\setminus\mathcal{A}^0}R(\alpha;\bolds\theta_\alpha) )$.

\begin{longlist}[(ii)]
\item[(i)] If $ |\mathcal{A}^0 |=0$, then $\hat\alpha_{\tau_n}$ is asymptotically loss efficient.
\item[(ii)] If $ |\mathcal{A}^0 |\geq 1$, then $\hat\alpha_{\tau_n}$ is consistent.
\end{longlist}
\end{thmm}

%re5.5 #&#
\begin{rem}
As in the one-dimensional case, the assumption,
$\hat{\bolds\theta}(\alpha)\mathop{\rightarrow}\limits^{P}\bolds\theta_\alpha$;
$\alpha\in\mathcal{A}\setminus\mathcal{A}^0$, is generally
satisfied. In fact, the assumption is guaranteed to hold when
$R(\alpha;\bolds\theta_0)=o(n^{(1+\delta)/2})$, for any
$\alpha\in\mathcal{A}$; see Lemma~A.5 of Chang, Huang and Ing (\citeyear{supp}).
\end{rem}

Here we consider only a multiplicative exponential model because of
two difficulties. First, for the two-dimensional exponential
covariance model, the asymptotic distribution of the ML estimate of
$(\sigma^2\kappa,\kappa)'$ is needed but has yet to be derived
unless $\kappa$ is assumed known [Du, Zhang and Mandrekar (\citeyear{Du2009}), \citet{Wang2011}]. Second, our proof relies on a decomposition of the
log-likelihood into different layers having different orders of
magnitude. Such a decomposition requires an innovative treatment of
the log-likelihood for the two-dimensional exponential model.
Further research is needed to characterize the asymptotic behavior
of GIC under the two-dimensional exponential covariance model or the
more general Mat\'{e}rn covariance model [Mat\'{e}rn (\citeyear{Matern1986})], but is
beyond the scope of this paper.
%s6 #&#
\section{Summary and discussion}
\label{discussion}

In this article, we study the asymptotic properties of GIC for
geostatistical model selection regardless of whether the covariance
model is correct or wrong, and establish conditions under which GIC
is consistent and asymptotically loss efficient. Some specific
examples that satisfy the regularity conditions are also provided.
To the best of our knowledge, this research is the first to provide
such results for GIC in geostatistical regression model selection.

The method we developed also sheds some light for solving linear
mixed model selection problems involving parameters that cannot be
estimated consistently. For example, consider a simple Laird--Ware
model [\citet{Laird1982}],
%
%e6.1 #&#
\begin{equation}
Z_{ij} = \mathbf{x}_{ij}'\bolds\beta +
\eta_i + \epsilon_{ij};\qquad i=1,\ldots,m, j=1,
\ldots,n_i , \label{Laird-Ware model}
\end{equation}
where $\mathbf{x}_{ij}$'s are $p$-vector of fixed effects, and
$\eta_i\sim N(0,\sigma^2)$ is the random effect for subject $i$,
independent of $\epsilon_{ij}\sim N(0,v^2)$. Here
$\bolds{\beta}\in\mathbb{R}^p$ is the regression-coefficient vector, and
$\bolds{\theta}=(\sigma^2,v^2)'$ consists of random-effect parameters.
Clearly, $\sigma^2$ in (\ref{Laird-Ware model}) cannot be estimated
consistently when $m$ is fixed [\citet{Longford2000}]. Nevertheless, as
shown below, it is still possible to derive a condition analogous to (C2). For simplicity, we
consider a simple case of (\ref{Laird-Ware model}) with mean zero
and no fixed effect. Let $\hat{\bolds{\theta}}$ be the ML estimate of
$\bolds{\theta}$ and $\bolds\theta_0=(v_0^2,\sigma_0^2)'$ be the true
parameter value. Applying an argument similar to that used to prove (2.10) of Chang, Huang and Ing (\citeyear{Chang2013}), twice the negative log-likelihood
of (\ref{Laird-Ware model}) can be written as
%
%e6.2 #&#
\begin{equation}
\qquad -2\ell(\bolds\theta) = n\log(2\pi) + \sum_{j=1}^m
\log n_j + n\log v^2 + n\frac{v_0^2}{v^2} + h(\bolds
\theta) + O_p(1), \label{log like fun LW model}
\end{equation}
where $  n = \sum_{i=1}^m n_i$, $
h(\bolds\theta)
=\sum_{i=1}^m \{\bolds\epsilon'_i\bolds\Sigma_i^{-1}\bolds\epsilon_i
-\mathrm{E} (\bolds\epsilon'_i\bolds\Sigma_i^{-1}\bolds\epsilon_i ) \}$,
$\bolds\epsilon_i=(\epsilon_{i1},\ldots,\epsilon_{i,n_i})'$ and
$\bolds\Sigma_j=\sigma^2\mathbf{1}_{n_j}\mathbf{1}_{n_j}' +v^2\mathbf{I}_{n_j}$. We
shall show that $\ell(\hat{\bolds\theta}) = \ell(\bolds\theta_0) +
O_p(1)$. Applying an argument similar to that used to prove Theorem~2.2 in Chang, Huang and Ing (\citeyear{Chang2013}),
%
%e6.3 #&#
\begin{equation}
\hat{\bolds\theta}=\bigl(v_0^2,\sigma_0^2
\bigr)'+ \bigl(O_p\bigl(n^{-1/2}
\bigr),O_p(1) \bigr)'. \label{Laird-Ware model:ML estimate}
\end{equation}

Let $\Theta_n=\{\bolds\theta\in\Theta\dvtx|\sigma^2-\sigma_0^2|<M,
|v^2-v_0^2|\leq M n^{-1/2} \}$ for any constant \mbox{$M>0$}. By Lemma B.1
of \citet{Chan2011} and an argument similar that used to prove~(2.12)
in Chang, Huang and Ing (\citeyear{Chang2013}), we have
\begin{eqnarray*}
&&\mathrm{E}\Bigl(\sup_{\bolds\theta\in\Theta_n}\bigl|h(\bolds\theta) - h(\bolds
\theta_0)\bigr|^2\Bigr)
\\
&&\qquad\leq  \sup_{\bolds\theta\in\Theta_n} \biggl\{\bigl(v^2-v_0^2
\bigr)^2\operatorname{var} \biggl(\frac{\partial}{\partial v^2}h(\bolds\theta)
\biggr) +\bigl(\sigma^2-\sigma_0^2
\bigr)^2\operatorname{var} \biggl(\frac{\partial}{\partial \sigma^2}h(\bolds\theta)
\biggr) \biggr\}
\\
&&\qquad= O(1),
\end{eqnarray*}
which implies $h(\hat{\bolds\theta})-h(\bolds\theta_0)=O_p(1)$. This together with
(\ref{log like fun LW model}) and (\ref{Laird-Ware model:ML
estimate}) gives $\ell(\hat{\bolds\theta}) = \ell(\bolds\theta_0) + O_p(1)$, indicating
some possibility to establish the asymptotic theory of GIC
for the Laird--Ware model, even when some random-effect parameter cannot be consistently estimated.

In this article, we focus only on variable selection under a certain covariance model. Clearly,
simultaneous selection of both variables and covariance models is an
interesting problem that deserves further investigation. Although we
believe that the framework we developed in this article can be
generalized to this problem, it will require introducing more
complex notation.

In addition, some more efforts are needed to completely characterize
GIC, even for the exponential covariance model in one dimension. We
note that both the candidate regressors in Examples
\ref{exm:white-noise} and \ref{exm:exp var} are not of bounded
variation (BV), whereas the polynomial regressors given by Example~\ref{exm:poly} are BV functions. It is of interest to know if BV
plays an important role. We conducted a small test simulation
experiment under the setup of (\ref{setup}) with only one regressor
$x(\cdot)$ and $v^2=0.5$, where $\mu(s)=1+x(s)$, $\eta(\cdot)$ is
given by (\ref{exp cov fun 1 dim}) with $\sigma^2=0.5$ and
$\kappa=1$, and data are sampled at $\{1/n,2/n,\ldots,1\}$. We
consider two functions for $x(\cdot)$, which are $f_1(s)
=s^2\sin(\pi/s)$ and $f_2(s) = s\sin(\pi/s)$, in combination with
three different sample sizes ($n=100, 500, 1000$). Note that
$f_1(\cdot)$ is of bounded variation on $[0,1]$, and $f_2(\cdot)$ is
not. The results based on 100 simulation replicates with known
$\sigma^2$, $\kappa$ and $v^2$ are shown in Table~\ref{spatial
dependent}. Clearly, GIC has better
ability in identifying the correct model when $f_2(\cdot)$, rather
than $f_1(\cdot)$, is used as the regressor, which partially supports
that BV may be an important factor.

%t1 #&#
\begin{table}
\caption{Frequencies of models selected by BIC
based on 100 simulation replicates,
where $\varnothing$ denotes the intercept only model and
$\alpha^0$ denotes the correct model}\label{spatial dependent}
\begin{tabular*}{\textwidth}{@{\extracolsep{4in minus 4in}}lcccc@{}}
\hline
 & \multicolumn{2}{c}{$\bolds{\mu(s)=1 + s^2\sin(\pi/s)}$} &
\multicolumn{2}{c@{}}{$\bolds{\mu(s)=1+ s\sin(\pi/s)}$}\\[-6pt]
 & \multicolumn{2}{c}{\hrulefill} &
\multicolumn{2}{c@{}}{\hrulefill}\\
$\bolds{n}$ & $\bm{\varnothing}$ & $\bolds{\alpha^0}$ &
$\bm{\varnothing}$ & \multicolumn{1}{c@{}}{$\bolds{\alpha^0}$} \\
\hline
 \phantom{0}100 & 67 & 33 & 38 & 62\\
 \phantom{0}500& 66 & 34 & 23 & 77\\
 1000& 76 & 24 & \phantom{0}8  & 92\\
\hline
\end{tabular*}
\end{table}


\begin{supplement}[id=suppA]
\stitle{Supplement to ``Asymptotic theory of generalized information
criterion for geostatistical regression model selection''\\}
\slink[doi]{10.1214/14-AOS1258SUPP}  %[doi,text={...}] - jei reikia suskaldyti doi
\sdatatype{.pdf}
\sfilename{aos1258\_supp.pdf}
\sdescription{The supplement materials contain the proofs of all
theorems in Section~\ref{section:examples}.}
\end{supplement}

% imsref loaded by akundreckaite, 2014-08-21 15:56:55

\printaddresses
\end{document}